\documentclass[12pt]{amsart}
\usepackage{amsmath}
\usepackage{amsfonts}
\usepackage{amsthm, upref}
\usepackage{graphicx}
 \usepackage{enumerate}
\usepackage[usenames, dvipsnames]{color}

    \numberwithin{equation}{section}

\usepackage{bm}
 \bmdefine\alphab{\mathbf{\alpha}}
\bmdefine\betab{\mathbf{\beta}}
\bmdefine\pib{\mathbf{\pi}}
\bmdefine\xib{\mathbf{\xi}}
\bmdefine\sigmab{\mathbf{\sigma}}

\newcommand{\comment}[1]{}
\newcommand{\eq}{\begin{equation}}
\newcommand{\en}{\end{equation}}
\newcommand{\rr}{\mathbb{R}}
\newcommand{\pp}{\mathbb{P}}
\newcommand{\nn}{\mathbb{N}}

\newcommand{\mcal}[1]{\mathcal{#1}}

\newcommand{\ev}{\mathbb E}

\newcommand{\ep}{\hfill $\Box$}

\newcommand{\NN}{\mathbb{N}}

\newcommand{\centS}{\overline{S}}
\newcommand{\centX}{\overline{X}}
\newcommand{\centx}{\overline{x}}
\newcommand{\mody}{\widetilde{\mathfrak{y}}}

\newcommand{\mR}{\mathfrak{R}}
\newcommand{\modR}{\widetilde{\mathfrak{R}}}
\newcommand{\mY}{\mathfrak{Y}}
\newcommand{\my}{\mathfrak{y}}
\newcommand{\modY}{\widetilde{\mathfrak{Y}}}
\newcommand{\mA}{\mathfrak{A}}

\newcommand{\wedgen}{\mathbb{H}}
\newcommand{\rwalk}{\Gamma}

\begin{document}

\theoremstyle{plain}
\newtheorem{thm}{Theorem}
\newtheorem{rem}{Remark}
\newtheorem{lemma}[thm]{Lemma}
\newtheorem{prop}[thm]{Proposition}
\newtheorem{cor}[thm]{Corollary}

\theoremstyle{definition}
\newtheorem{defn}{Definition}
\newtheorem{cond}{Condition}
\newtheorem{asmp}{Assumption}
\newtheorem{notn}{Notation}
\newtheorem{prb}{Problem}

\theoremstyle{remark}
\newtheorem{rmk}{Remark}
\newtheorem{exm}{Example}
\newtheorem{clm}{Claim}

\title[Asymmetric collisions]{Systems of Brownian particles with asymmetric collisions}

\author{Ioannis Karatzas, Soumik Pal, and Mykhaylo Shkolnikov}
\address{INTECH Investment Management \\ One Palmer Square\\ Princeton, NJ 08542 and Columbia University \\ Department of Mathematics\\ New York, NY 10027}
\email{ik@enhanced.com, ik@math.columbia.edu}
\address{Department of Mathematics\\ University of Washington\\ Seattle WA 98195}
\email{soumikpal@gmail.com}
\address{INTECH Investment Management \\ One Palmer Square\\ Princeton, NJ 08542 and University of California, Berkeley\\ Department of Statistics\\ Berkeley, CA 94720}
\email{mshkolni@gmail.com}

\thanks{Research   partially supported by NSF grants DMS-09-05754 (I. Karatzas) and DMS-10-07563 (S. Pal).}

\keywords{Determinantal processes, Interacting particle systems, Invariance principles, Reflected Brownian motions, Skorokhod maps, Stochastic Portfolio Theory, Strong solutions of stochastic differential equations, Triple collisions}

\subjclass[2000]{60K35, 60H10, 91B26}

\date{\today}

\begin{abstract}
We study systems of Brownian particles on the real line, which interact by splitting the local times of collisions among themselves in an asymmetric manner. We prove the strong existence and uniqueness of such processes and identify them with the collections of ordered processes in a Brownian particle system, in which the drift co\"efficients, the diffusion co\"efficients, and the collision local times for the individual particles are assigned according to their ranks. These 
%Soumik 
Brownian systems 
can be viewed as generalizations of those arising in first-order models for equity markets   in the context of stochastic portfolio theory, and are able to correct for several shortcomings of such models  while being equally amenable to computations. We also show that, in addition to being of interest in their own right, such systems of Brownian particles arise as universal scaling limits of systems of 
%random walks 
% Yannis
jump processes on the integer lattice with local interactions. 
%Soumik
A key step in the proof is the analysis of a generalization of Skorokhod maps which include `local times' at the intersection of faces of the nonnegative orthant.  The result extends the convergence of TASEP to its continuous analogue. 
Finally, we identify those among the Brownian particle systems which have a probabilistic structure of determinantal type.
\end{abstract}

\maketitle

%%%%%%%%%%%%%%%%%%%%%%%%%%%%%%%%%%%%%%%
\section{Introduction}
\label{sec1}
%%%%%%%%%%%%%%%%%%%%%%%%%%%%%%%%%%%%%%%

Systems of Brownian particles with various types of interactions have
%Soumik
been widely studied. More recently, 
Brownian particles with electrostatic repulsion (\textsc{Dyson}'s Brownian motion, see for instance  section 4.3 in \cite{AGZ}) play a central r\^ole in the understanding of the universality properties of large Hermitian random matrices with independent entries (\textsc{Wigner} matrices, see the survey \cite{EY} and the references there). In addition, systems of Brownian particles interacting through their ranks, which were originally introduced in the context of the piecewise-linear filtering problem and in the resulting study of diffusions with piecewise-constant characteristics (see \cite{BP}), have been of great importance in the study of large equity markets within stochastic portfolio theory (see \cite{F}, \cite{FK}). Finally, systems of Brownian particles interacting by stickiness have been recently introduced as continuous analogues of a certain random evolution for the distribution of mass on the integer lattice (see \cite{HW} and the references there).   

\medskip
As a starting point towards the formulation of our setup, we recall from section 13.1 in \cite{FK} that the ordered processes in a system of Brownian particles interacting through their ranks are given by independent Brownian motions, which have constant drift and diffusion co\"efficients and collide in a \textit{symmetric} fashion; that is,  each collision local time is   split equally between the two colliding particles. In contrast, the particles in the process introduced by \textsc{Warren} in \cite{Wa} evolve as independent Brownian motions which have constant drift and diffusion co\"efficients and collide in a \textit{totally asymmetric} manner; that is, the collision local time is assigned  entirely to one of the two colliding particles. 

%\medskip

%%%%%%%%%%%%%%%%%%%%%%%%%%%%%%%%%%%%%%%
\subsection{The Setup}
\label{sec0.0}
%%%%%%%%%%%%%%%%%%%%%%%%%%%%%%%%%%%%%%%

%Yannis: I thought it was a good idea to separate this discussion of the "setup" as a new subsection. 
The  dynamics of (\ref{ranks}) below  give  the formal description of a system of 
%Soumik
ordered  
Brownian particles
%Soumik
on the line, 
 which move as independent Brownian motions with constant drift and diffusion co\"efficients and collide \textit{asymmetrically}; that is, the collision local times are apportioned unequally, in a manner that depends on the ranks of the particles involved in the collisions. 

Consider a continuous, $\,n-$dimensional  semimartingale $\,(R_1(\cdot),R_2(\cdot),\ldots, R_n(\cdot))\,$ with values in the Weyl chamber $\, \mathbb{W}^n = \big\{ (r_1,r_2,\ldots, r_n): \infty > r_1 \ge r_2 \ge \cdots \ge r_n > - \infty \big\}\,$, and with dynamics of the form 
%the following SDE for the ranks
\eq
\label{ranks}
%\mathrm{d}
~R_k(t)\,=\, R_k(0) + b_k\,  %\mathrm{d}
t+\sigma_k\,  %\mathrm{d}
\betab_k(t)+q_k^-\,  %\mathrm{d}
\Lambda^{(k,k+1)}(t)-q_k^+\,  %\mathrm{d}
\Lambda^{(k-1,k)}(t)\,,\quad 0 \le t < \infty
\en
for $\,k=1,2,\ldots,n\,$. Here the drifts $\,b_1,\,b_2,\ldots,\,b_n\,$ are given real numbers; the dispersions $\,\sigma_1,\,\sigma_2,\ldots,\,\sigma_n\,$ are given  positive real numbers;   the collision parameters $\,q_1^{\pm},\,q_2^{\pm},\ldots,q_n^{\pm}\,$ are given positive real numbers satisfying 
\eq
\label{elastic}
q_k^- + q_{k+1}^+\,=\,1,\quad k=1,2,\ldots,n-1\,;
\en
and the processes $\,\betab_1 (\cdot) ,\, \betab_2 (\cdot),\ldots,\,\betab_n(\cdot)\, $ are independent standard Brownian motions. On the other hand, for each $\,k=1,2,\ldots,n-1\,$ the process 
\eq\label{eq:whatislambda}
\Lambda^{(k,k+1)} (\cdot) \,\equiv \, L^{R_k   - R_{k+1}} (\,\cdot\,; 0)
\en
is the right-sided local time accumulated at the origin  by the nonnegative semimartingale $\,R_k (\cdot) - R_{k+1} (\cdot) \,$; we set $\,\Lambda^{(0,1)} (\cdot)\equiv\Lambda^{(n,n+1)} (\cdot) \equiv0\,$. The ``regulating" r\^ole of these local times in (\ref{ranks}) is to make sure the  resulting    process $\,   ( R_1(\cdot), R_2(\cdot),\ldots, R_n(\cdot) )=:R(\cdot) \,$ takes values in the wedge $\,\mathbb{W}^n\,$ at all times. This process $\, R(\cdot)\,$ can thus be regarded as  Brownian motion   with reflection on the faces of the polyhedral domain $\, \mathbb{W}^n\,$, in the sense of  \textsc{Harrison \& Williams} (1987) (cf.$\,$\cite{HW1}, \cite{W}). 

\smallskip

As we show in section 3, Brownian particle systems of the type \eqref{ranks} are not only of interest in their own right, but also arise as universal scaling limits for systems of %random walks 
jump processes 
%Soumik
%\footnote{~ I changed ``random walk" here to ``jump process" (IK).} 
%
on the integer lattice with local interactions. Consider a system of $n$ particles on the integer lattice, moving according to (possibly asymmetric) continuous time simple %random walks
jump processes, which are independent as long as the particles are located at $n$ different sites. When two or more particles land at the same site (we will refer to such events as ``collisions''),   the jump rates of the particles change, and   in a manner  that preserves the order of the particles. As we explain in section 3, such particle systems converge under a diffusive rescaling of time and space to the solution of %an SDE 
a stochastic equation of the type \eqref{ranks}. The proof of this result necessitates a detailed study of 
%Soumik
Skorokhod maps that transform noise to processes constrained to stay in the nonnegative orthant, but might involve `local time' push from the intersection of multiple faces of the orthant. 
Along the way, we generalize the boundary property of reflected Brownian motion established in \cite{RW} (see Lemma \ref{genRW} below) and the invariance principle for reflected Brownian motion of \cite{W1} (see Proposition \ref{inv_prin} below).        
 
\medskip
We now discuss how the processes in \eqref{ranks} can be seen as  describing the order statistics in Brownian particle systems, in which the particles are allowed to exchange their ranks. As we explain below, the latter can be used as models for the logarithmic capitalizations in large equity markets, and generalize the so-called ``first order models" of stochastic portfolio theory. Consider an $\,n-$dimensional process $\,(X_1(\cdot),X_2(\cdot),\ldots,X_n(\cdot))\,$ that satisfies the system of stochastic differential equations  
\eq
\label{names}
\begin{split}
\mathrm{d}X_i(t)\,= \, \sum_{k=1}^n \,\mathbf{1}_{\{X_i(t)=R^X_k(t)\}}\,b_k\,\mathrm{d}t 
+\sum_{k=1}^n \mathbf{1}_{\{X_i(t)=R^X_k(t)\}}\,\sigma_k\,\mathrm{d}W_i(t) \\
+\sum_{k=1}^n \mathbf{1}_{\{X_i(t)=R^X_k(t) \}}\, \big( q_k^--(1/2) \big)\,\mathrm{d}\Lambda^{(k,k+1)}(t)\\
-\sum_{k=1}^n \mathbf{1}_{\{X_i(t)=R^X_k(t) \}}\, \big(q_{k }^+ -(1/2)\big)\,\mathrm{d}\Lambda^{(k-1,k)}(t)\,.
\end{split}
\en
Here $\, W_1(\cdot),W_2(\cdot),\ldots,W_n(\cdot)\, $ are independent standard Brownian motions; the ``descending order statistics" 
\eq
\label{os}
 \max_{1 \le j \le n } X_j (\cdot) =: R^X_1 (\cdot) \,\ge \, R^X_2 (\cdot) \, \ge\,  \ldots   \, \ge \, R^X_n (\cdot) := \min_{1 \le j \le n } X_j (\cdot) \,
 \en
  are the ranked processes corresponding to the particle system $\,X_1 (\cdot),X_2(\cdot),\ldots,X_n(\cdot)\,$ with lexicographic resolution of ties; whereas  
  \eq
  \label{CollLocTime}
  \Lambda^{(k,\ell)} (\cdot) \,\equiv \, L^{R^X_k   - R^X_{\ell}} (\cdot\,; 0) \,, \qquad \ell \ge k+1
\en
denotes the local time accumulated at the origin  by the nonnegative semimartingale $\,R^X_k (\cdot) - R^X_{\ell} (\cdot) \,$.
  
\smallskip
  
Assume now that a weak solution to the system of \eqref{names} exists and satisfies 
\eq\label{notsticky}
{\mathcal Leb}\,\big(\{t\geq0\,|\;\exists \,\,1\leq i<j\leq n:\;X_i(t)=X_j(t)\}\big)\,=\,0 
\en
almost surely, where $\,{\mathcal Leb}\,$ denotes the Lebesgue measure on $[0,\infty)$. Then with the notation $\,\mathcal{N}_k(t)=|\{ i: X_i (t) = R^X_k (t)\}|\,$ for the number of particles occupying the $\,k$-th rank at time $\,t\,$, the \textsc{Banner \& Ghomrasni} (2008) formula    
\eq
\label{BaG}
\mathrm{d}R^X_k(t)\,=\, \sum_{i=1}^n \mathbf{1}_{\{R^X_k(t)=X_i(t)\}}\,\mathrm{d}X_i(t)+ \frac{1}{\,\mathcal{N}_k (t) \,}\,\left( \, \sum_{\ell = k+1}^{n} \mathrm{d}\Lambda^{(k,\ell)}(t)- \sum_{\ell =1}^{k-1}\mathrm{d}\Lambda^{(\ell,k)}(t) \right)   
\en
(cf.$\,$Theorem 2.3 in \cite{BG}) shows that the process of spacings 
\eq
\label{spac}
Z^X (\cdot)\,:=\,\big(R^X_1 (\cdot)-R^X_2(\cdot)\,,\,R^X_2 (\cdot)-R^X_3 (\cdot)\,,\cdots,\,R^X_{n-1} (\cdot)-R^X_n (\cdot)\big), 
\en 
when away from the boundary of the nonegative orthant $(\rr_+)^{n-1}$, moves according to the 
%Soumik
multidimensional process
 $\,((b_1-b_2)t+\sigma_1\,\betab^X_1(t)-\sigma_2\,\betab^X_2(t),\ldots,(b_{n-1}-b_n)t+\sigma_{n-1}\,\betab^X_{n-1}(t)-\sigma_n\,\betab^X_n(t))\,$, $t\geq0$. Here the processes 
\eq
\label{beta}
\betab^X_k (\cdot) \,:=\,  \sum_{i=1}^n \int_0^{\,\cdot} \mathbf{1}_{\{X_i(t)=R^X_k(t)\}}\, \mathrm{d}W_i(t)\,,\qquad k=1,2,\ldots,n 
\en
are independent standard Brownian motions, by virtue of the \textsc{P. L\'evy} theorem (see section 3 in \cite{BFK} for a very similar derivation). But on the strength of   Lemma \ref{genRW} below, which generalizes the boundary property of reflected Brownian motion established by \textsc{Reiman \& Williams} (1988) in \cite{RW}, the triple- or higher-order collision local times $\,\Lambda^{(k,\ell)}(\cdot)\,$ vanish for all $\,\ell \ge k+2\,$, so the expression in the \textsc{Banner-Ghomrasni} (2008) formula (\ref{BaG}) simplifies to 
\begin{eqnarray*}
\mathrm{d}R^X_k(t)&=&\sum_{i=1}^n \mathbf{1}_{\{R^X_k(t)=X_i(t)\}}\,\mathrm{d}X_i(t)+ \frac{1}{\,2\,}\,\mathrm{d}\Lambda^{(k,k+1)}(t)-\frac{1}{\,2\,}\,\mathrm{d}\Lambda^{(k-1,k)}(t)   \\
&=& \sum_{i=1}^n \mathbf{1}_{\{R^X_k(t)=X_i(t)\}}\,\Big( \sum_{\ell=1}^n \,\mathbf{1}_{\{X_i(t)=R^X_\ell(t)\}}\,b_\ell\,\mathrm{d}t 
+\sum_{\ell=1}^n \mathbf{1}_{\{X_i(t)=R^X_\ell(t)\}}\,\sigma_\ell\,\mathrm{d}W_i(t)\Big) \\
&&+ \sum_{i=1}^n \mathbf{1}_{\{R^X_k(t)=X_i(t)\}}\,\sum_{   \ell =1}^n 
\mathbf{1}_{\{X_i(t)=R^X_\ell (t) \}}\, \big( q_\ell^-- (1/2) \big)\,\mathrm{d}\Lambda^{(\ell,\ell+1)}(t) \\
&&- \sum_{i=1}^n \mathbf{1}_{\{R^X_k(t)=X_i(t)\}}\,\sum_{ \ell=1}^n
\mathbf{1}_{\{X_i(t)=R^X_\ell(t) \}}\, \big( q_{\ell }^+ -(1/2)\big)\,\mathrm{d}\Lambda^{(\ell-1,\ell)}(t) \\
&&+ \frac{1}{\,2\,}\,
\mathrm{d}\Lambda^{(k,k+1)}(t)-\frac{1}{\,2\,}\,\mathrm{d}\Lambda^{(k-1,k)}(t)\,,\quad k=1,\ldots,n\,. 
\end{eqnarray*}
Evaluating the sums of this expression over $i$ and $\ell$ and recalling the notation of (\ref{beta}), we get the dynamics \eqref{ranks} for the descending order statistics of (\ref{os}). 

%%%%%%%%%%%%%%%%%%%%%%%%%%%%%%%%%%%%%%%
\subsection{Interpretation}
\label{sec1.0}
%%%%%%%%%%%%%%%%%%%%%%%%%%%%%%%%%%%%%%%

We shall think of the processes $\,X_1(\cdot),X_2 (\cdot),\ldots,X_n(\cdot)\,$ as representing the positions of a finite collection of Brownian particles, whose    drift  and dispersion co\"efficients are assigned according to the ranks   occupied by the particles when ordered as in (\ref{os}). When the particles collide, they interact  asymmetrically with their nearest neighbors through the collision local times at the origin, in the specific manner of (\ref{names}) and with the notation of   (\ref{CollLocTime}). We shall speak of the $\,\mathbb{R}^n-$valued semimartingale $\,(X_1(\cdot),\ldots,X_n(\cdot))\,$ as the process of ``names" (positions of individual particles), and of the components of the $\,\mathbb{W}^n-$valued semimartingale $\,(R_1(\cdot), \ldots,R_n(\cdot))\,$ as the associated ``ranked processes"   (descending order statistics). 

\smallskip

If we denote by $\,\mathfrak{r}_t (i)\,$ the rank occupied by particle $\,i\,$ at time $\,t\,$, we can write the system of equations (\ref{names}) in the very informal, yet suggestive and slightly more compact form 
\eq
\begin{split} \label{names too}
\mathrm{d}X_i(t)\,=\,b_{\,\mathfrak{r}_t (i)}\,\mathrm{d}t  
+\sigma_{\mathfrak{r}_t (i)}\,\mathrm{d}W_i(t) 
+\big(q_{\mathfrak{r}_t (i)}^- - (1/2)\big)\,\mathrm{d}\Lambda^{ (\mathfrak{r}_t (i), \,\mathfrak{r}_t (i)+1)}(t)\\
~~-\,\big( q_{\mathfrak{r}_t (i)}^+ - (1/2)\big)\,\mathrm{d}\Lambda^{(\mathfrak{r}_t (i)-1,\, \mathfrak{r}_t (i))}(t)\,.
\end{split}
\en
To wit: at any given time $\,t\,$, every particle $\,i\,$ gets assigned   drift and   dispersion parameters according to its current rank $\, \mathfrak{r}_t (i)  \,$, and feels an upward (respectively, downward) local-time-like pressure, or ``drag",  when colliding with  the particle right below it (respectively, right above  it)  in proportion to $\,   q^{\mp}_{\mathfrak{r}_t (i)}-(1/2)\,$. %, with $\,k=\mathfrak{r}_t (i)\,$ being its current rank. 

%%%%%%%%%%%%%%%%%%%%%%%%%%%%%%%%%%%%%%%
\subsection{Some Special Cases}
\label{sec1.1}
%%%%%%%%%%%%%%%%%%%%%%%%%%%%%%%%%%%%%%%
 
It is instructive to compare the system \eqref{names} in the two-dimensional case $\,n=2\,$ with the systems of equations studied by \textsc{Fernholz, Ichiba \& Karatzas} (2012) in \cite{FIK} (see the systems of equations (1.2), (1.3) and (4.13), (4.14) of that paper). As one can see by comparing the co\"efficients of the local time terms, in the case $\,n=2\,$  the system \eqref{names} is a special case of the system of equations (4.13)--(4.14) in \cite{FIK}, namely 
$$
\mathrm{d} X_1 (t) = \mathbf{ 1}_{ \{ X_1 (t) \ge X_2 (t) \}} \big( b_1   \mathrm{d}t + \sigma_1  \mathrm{d} W_1 (t) \big)+ \mathbf{ 1}_{ \{ X_1 (t) < X_2 (t) \}} \big( b_2   \mathrm{d}t + \sigma_2  \mathrm{d} W_2 (t) \big)+ \kappa \,  \mathrm{d} L ^{|X_1 - X_2|} (t)
$$
$$
\mathrm{d} X_2 (t) = \mathbf{ 1}_{ \{ X_1 (t) < X_2 (t) \}} \big( b_1   \mathrm{d}t + \sigma_1  \mathrm{d} W_1 (t) \big)+ \mathbf{ 1}_{ \{ X_1 (t) \ge  X_2 (t) \}} \big( b_2   \mathrm{d}t + \sigma_2  \mathrm{d} W_2 (t) \big)+ \kappa \,  \mathrm{d}  L ^{|X_1 - X_2|} (t)
$$

\medskip
\noindent
with $\,\kappa = q^-_1-(1/2)=(1/2)-q^+_2\,$. The inequalities in the indicators reflect the convention on lexicographic resolution of ties we referred to earlier. 

\smallskip

In particular, with  $\, \Upsilon (\cdot) = X_1 (\cdot) - X_2 (\cdot)\,$, $\, \Xi (\cdot) = X_1 (\cdot) + X_2 (\cdot)\,$, $\, \lambda_1 = b_1 - b_2\,$ and $\,\lambda_2=b_1+b_2\,$, we have $\, R_1 (\cdot) - R_2 (\cdot) = | \Upsilon (\cdot)|\,$, and the processes 
$$
W( t) = \Upsilon ( t) + \lambda_1 \int_0^{ t} \mathrm{sgn} (\Upsilon(s))\, \mathrm{d} s\,, \quad V( t) = \Xi( t) - \lambda_2 \, t -2\, \kappa\, L ^{|\Upsilon|} (t)\,, \quad 0 \le t < \infty
$$
are now Brownian motions with diffusion co\"efficients $\, \sqrt{\sigma_1^2+\sigma_2^2\,}\,$ and with covariation $\, \langle W, V \rangle (\cdot) = (\sigma_1^2 - \sigma_2^2)   \int_0^{ \, \cdot} \mathrm{sgn} (\Upsilon(t))\, \mathrm{d} t\,$. 

\smallskip

Let us also remark that the solution to the above two-dimensional system can  be realized as the solution of the system of equations (1.2)--(1.3) in \cite{FIK} with 
\eq
\kappa \,= \,q_1^--\frac{ 1}{\,2\,}\,=\,\frac{1}{\,2\,}-q_2^+
\,=\,\frac{1-\eta_1}{2}\,=\,\frac{1-\eta_2}{2}\,=\,\frac{1-\zeta_1}{2}\,=\,\frac{1-\zeta_2}{2}\,,
\en  
and with $\eta_1,\,\eta_2,\,\zeta_1,\,\zeta_2$ the parameters introduced in \cite{FIK}. In particular, the condition (1.5) in \cite{FIK}, which is necessary and sufficient for the well-posedness of the system of equations (1.2)--(1.3) in \cite{FIK}, is here fulfilled.  

\medskip
\noindent $\bullet~$ Now, suppose that we have $$b_1=\ldots=b_n=0\, , ~~\sigma_1=\ldots=\sigma_n=1\,~~~~\; \hbox{and} ~\;~\,q^-_1=\ldots=q^-_n=0\,, ~~\,q^+_1=\ldots=q^+_n=1\,.$$ In this case one can think of the particles as having masses that decrease  from right to left, so that the mass of each particle is negligible compared to the mass of its right neighbor. To wit, whenever a particle collides with its right neighbor, it is reflected off this considerably ``heavier" particle. In this situation, the vector of ranked processes $\,R(\cdot)\,$ is given by the continuous version of the {\it Totally Asymmetric Simple Exclusion Process} (TASEP). We refer the reader to section 4 in \cite{Wa} for some of the properties of this process, and to \cite{GS} for its appearance as the scaling limit of TASEP.   

\medskip
\noindent $\bullet~$ Let us also note that the special case
  \eq
\label{half}
q^{\,\pm}_k\,=\, \frac{1}{\,2\,}\,, \qquad k=1,2,\ldots,n\,,
 \en
in which all local times disappear from (\ref{names}),  (\ref{names too}) and  get equal weights in (\ref{ranks}), and {\it all collisions of  ranked   particles are symmetric,} was studied in detail by \textsc{Ichiba, Karatzas \& Shkolnikov} in \cite{IKS}; in this case, {\it individual}  particles collide with each other without feeling any local time drag from their nearest neighbors. 

%%%%%%%%%%%%%%%%%%%%%%%%%%%%%%%%%%%%%%%
\subsection{Outline}
\label{sec1.2}
%%%%%%%%%%%%%%%%%%%%%%%%%%%%%%%%%%%%%%%

The rest of the paper is organized as follows. In section 2.1 we show the strong existence and uniqueness of the solution to the system of equations \eqref{ranks}. Subsequently, in section 2.2 we prove the appropriate generalization of the boundary property of reflected Brownian motion established in \cite{RW} to processes involving local times accumulated on lower-dimensional faces of the boundary; then, we study the attainability of lower-dimensional parts of the boundary by such processes in section 2.3. These results are combined in section 2.4 to prove the strong existence and uniqueness of the solution to the system of equations \eqref{names}, under the explicit conditions on non-attainability of lower-dimensional boundaries derived in section 2.3. The latter results extend those of \cite{BP}, \cite{IKS} and \cite{SV} (within our setup, the first two only address the case \eqref{half}, whereas the last one only addresses the case of rank-independent parameters). Then, in section 2.5 we view the solution to \eqref{names} as a model for the logarithmic capitalizations in large equity markets and analyze the resulting capital distributions. 

\smallskip

Section 3 is devoted to the identification of the processes in \eqref{ranks} as universal scaling limits of systems of one-dimensional %random walks 
continuous time jump processes with local interactions. It includes an extension of the invariance principle of \cite{W1} in section 3.3 (see Proposition \ref{inv_prin}). Finally, in section 4 we characterize the sets of collision parameters, for which the solution to \eqref{ranks} has a probabilistic structure of determinantal type. These results generalize those in \cite{Wa}, on Brownian particle systems with totally asymmetric collisions.  

%%%%%%
\section{Analysis of the continuous process}
%Section added by Soumik 09/14

%%%%%%%%%%%%%%%%%%%%%%%%%%%%%%%%%%%%%%%
\subsection{Ranks}
%%%%%%%%%%%%%%%%%%%%%%%%%%%%%%%%%%%%%%%

We start with the construction of the vector $\,R(\cdot)=(R_1(\cdot),\ldots,R_n(\cdot))\,$ of ranked processes. We note first that, due to the positivity of the co\"efficients $q^{\pm}_1,q^{\pm}_2,\ldots,q^{\pm}_n$, there exist positive constants $c_1,c_2,\ldots,c_n$ such that   $\sum_{i=1}^n c_i R_i(t)$, $t\geq0$ is a Brownian motion with drift, that is, the contribution of the local times to its dynamics vanishes. This observation allows us to construct the ranked processes $\, R_1(\cdot),R_2(\cdot),\cdots,R_n(\cdot)\,$ using the following procedure: first, we define the auxiliary Brownian motion $\,\widetilde{R}(\cdot)=(\widetilde{R}_1(\cdot),\widetilde{R}_2(\cdot),\ldots,\widetilde{R}_n(\cdot))\,$, for which
$$
 \mathrm{d}\left(\sum_{k=1}^n c_k\widetilde{R}_k(t)\right)\,=\,\left(\sum_{k=1}^n c_k\,b_k\right)\mathrm{d}t \,+ \, \sum_{k=1}^n c_k\,\sigma_k\,\mathrm{d}\betab_k(t) 
 $$
 and
 $$
\mathrm{d}\big(\widetilde{R}_k(t)-\widetilde{R}_{k+1}(t)\big)\,=\,\big(b_k-b_{k+1}\big)\,\mathrm{d}t 
%&&\quad\quad\quad\quad\quad\quad\quad
+\,\sigma_k\,\mathrm{d}\betab_k(t)-\,\sigma_{k+1}\,\mathrm{d}\betab_{k+1}(t) %\,,~~\;k=1,2,\ldots,n-1\,.
$$

\medskip
\noindent
for $\,k=1,2,\ldots,n-1\,$. Next, we introduce   the process $$ Y(\cdot)\,:=\, \Big(\sum_{k=1}^n c_k\widetilde{R}_k (\cdot)\, ,\,\widetilde{R}_1 (\cdot) -\widetilde{R}_2(\cdot)\, ,\ldots,\,\widetilde{R}_{n-1} (\cdot)-\widetilde{R}_n(\cdot) \Big)$$   and apply the \textsc{Harrison-Reiman} (1981) version of the \textsc{Skorokhod} reflection map $\Psi^{HR}$ for the   orthant $\, ( \mathbb{R}_+)^{n-1}\,$ to the last $(n-1)$ components of $\,Y(\cdot)\,$, using the reflection matrix
\eq\label{refl_matrix}
\mathcal R=\mathbf{ I}_{\,n-1}-\mathcal Q\,, \quad \mathrm{where} \quad \mathcal Q\,:=\, 
\left(\begin{array}{llll}
0 & q_2^- & 0 & 0  \\
q_2^+ & 0 & q_3^- & 0  \\
0 & \ddots & \ddots & \ddots  \\
0 & 0  & q_{n-1}^+& 0
\end{array}\right) 
\en 
and $\, \mathbf{ I}_{\,n-1}\,$ is the unit $\, (n-1) \times (n-1)\,$ matrix. 

The main observation here, is that {\it the spectral radius of the matrix $\,{\mathcal Q}\,$ is strictly less than $1$}. Indeed, the transpose $\mathcal Q^T$ is an irreducible substochastic matrix, and can be made into a stochastic matrix by adding an absorbing point. Hence, by virtue of the \textsc{Perron-Frobenius} Theorem, the spectral radius of $\,{\mathcal Q}^T\,$, and hence also of $\,\mathcal Q\,$, is strictly less than $1$.

\smallskip

All in all, we see that Theorem 1 of \cite{HR} is applicable here. We can now complete the definition of the process $\, R(\cdot) = ( R_1(\cdot), \cdots, R_n(\cdot) )\,$ by imposing
$$
 \sum_{k=1}^n c_k\,R_k (\cdot)\, =\, \sum_{k=1}^n c_k\,\widetilde{R}_k (\cdot)\, ,
 $$
$$
\big(R_1 (\cdot)-R_2 (\cdot),\ldots,R_{n-1} (\cdot)-R_n (\cdot) \big)\,=\,\Psi^{HR}\big(\widetilde{R}_1 (\cdot)-\widetilde{R}_2 (\cdot),\ldots,\widetilde{R}_{n-1} (\cdot)-\widetilde{R}_n (\cdot)\big)\,. 
$$
The strong uniqueness of the so-contructed process $\,R (\cdot)\,$ follows   from the uniqueness of the solution to the multi-dimensional \textsc{Skorokhod} reflection problem  in \textsc{Harrison \& Reiman} (1981), \textsc{Reiman} (1984); see Theorem 1 of \cite{HR}. 

\medskip

Let us close this subsection by observing from (\ref{ranks})  that the process of spacings 
\eq\label{spacings}
Z  (\cdot)\,:=\, 
\big(R_1 (\cdot)-R_2(\cdot),\,R_2 (\cdot)-R_3 (\cdot),\cdots,\,R_{n-1} (\cdot)-R_n (\cdot)\big)
\en  
for the process $\, R(\cdot)\,$ we just constructed, is a reflected Brownian motion (RBM)  in the nonnegative orthant $\,(\rr_+)^{n-1}\,$,  with drift vector $\, (b_1 - b_2, \cdots, b_{n-1}-b_n)\,$, covariance matrix 
\eq\label{cov_matrix}
\mathcal{A}\,=\,\left(
\begin{array}{llll}
\sigma_1^2+\sigma_2^2 & -\sigma_2^2 & 0 & 0 \\
-\sigma_2^2 & \sigma_2^2+\sigma_3^2 & -\sigma_3^2 & 0 \\
0 & \ddots & \ddots & \ddots \\
0 & 0 & -\sigma_{n-1}^2 & \sigma_{n-1}^2+\sigma_n^2
\end{array}
\right),
\en
and reflection matrix $\,\mathcal R=\mathbf{ I}_{n-1}-\mathcal Q\,$ given as in \eqref{refl_matrix}. 

%%%%%%%%%%%%%%%%%%%%%%%%%%%%%%%%%%%%%%%
\subsection{A boundary property of reflected Brownian motion}
%%%%%%%%%%%%%%%%%%%%%%%%%%%%%%%%%%%%%%%

Throughout the paper, we rely many times on a generalization of the boundary property of reflected Brownian motion established in \cite{RW}, which is of interest in its own right. 

Consider a continuous semimartingale $Q(\cdot)$ taking values in the orthant $(\rr_+)^{n-1}$ and satisfying 
\eq\label{genRBM}
Q(\cdot)=B(\cdot)+\sum_{k=1}^{n-1} \mR^{(k)}\,\mY^{(k)}(\cdot),
\en
where $B(\cdot)$ is an $(n-1)$-dimensional Brownian motion with a constant drift vector and a constant, nondegenerate diffusion matrix. Here, for each $k\in\{1,2,\ldots,n-1\}$ and with $$\, m \,=\, \binom{n-1}{k}\,,$$ $\mR^{(k)}$ is an $\,(n-1)\times  m\,$ matrix with real entries; whereas  %for each $k\in\{1,2,\ldots,n-1\}$, 
$\mY^{(k)}(\cdot)$ is a continuous $(\rr_+)^{m}-$valued process, whose components are indexed by the sets $J\subset\{1,2,\ldots,n-1\}$ with $k$ elements, start at $\,\mY^{(k)}(0)=0\,$, are non-decreasing and satisfy
\begin{eqnarray}
\int_0^\infty \sum_{j\in J} \mathbf{1}_{\{Q_j(s)>0\}}\,\mathrm{d}\mY^{(k)}_J(s)=0,\quad J\subset\{1,2,\ldots,n-1\},\;\;|J|=k, \label{loccond1} \\
\begin{split}
\forall\,\,\,0\leq s<t:\;\mY^{(|J_2|)}_{J_2}(t)-\mY^{(|J_2|)}_{J_2}(s)\leq \mY^{(|J_1|)}_{J_1}(t)-\mY^{(|J_1|)}_{J_1}(s), \quad\quad\quad \label{loccond2}\\\hbox{for any}~~~
J_1\subset J_2\subset\{1,2,\ldots,n-1\}. 
\end{split}
\end{eqnarray}    
As we show in the following lemma, if the matrix $\,\mR^{(1)}\,$ is completely$-\mathcal{S}$ (see \cite{RW} for the definition and a characterization of completely$-\mathcal{S}$ matrices), then $Q(\cdot)$ can be identified with a reflected Brownian motion in the sense of \cite{HW1}, \cite{W}.

\begin{lemma}\label{genRW}
Let $Q(\cdot)$ be a process as in \eqref{genRBM}, and suppose that the matrix $\mR^{(1)}$ is completely$-\mathcal{S}$. Then all processes $\,\mY^{(k)}(\cdot)\,$, $k=2,3,\ldots,n-1$ are identically equal to zero.   
In particular, $Q(\cdot)$ is a reflected Brownian motion in the orthant $(\rr_+)^{n-1}$ with reflection matrix $\,\mR^{(1)}$. 
\end{lemma}

\noindent{\bf Proof:} The proof is similar to the proof of Theorem 1 in \cite{RW}. In view of Girsanov's Theorem, it suffices to consider the case that the drift vector of $B(\cdot)$ is equal to zero. We consider the functions $\,\phi_\varepsilon\ (\cdot)\,$, $  \varepsilon  \in(0,1)$ as in the proof of Lemma 4 in \cite{RW}; these functions are harmonic for  the generator of $B(\cdot)$ and bounded on compact subsets of $(\rr_+)^{n-1}$, uniformly over $ \,\varepsilon  \in(0,1)\,$. We claim that, for any $k=1,2,\ldots,n$ and any column $v$ of the matrix $\mR^{(k)}$, there is a constant $C_k<\infty$  depending only on the matrices $\,\mR\,$, $\mR^{(k)}\,$ and the diffusion matrix of $B(\cdot)$, such that
\eq\label{grad_bound}
\forall\,\varepsilon>0,\;x\in(\rr_+)^{n-1}:\quad v\cdot\nabla\phi_\varepsilon(x)\geq -C_k\,.
\en
Indeed, one can argue exactly as on page 94 of \cite{RW} in the derivation of the bound (24) there. Now, define the stopping times
\[
\tau_m\,=\,\min\Big(m\,,\,\inf\Big\{t\geq0:\;\|Q(t)\|+\sum_{k=1}^{n-1} \|\mY^{(k)}(t)\|\geq m\Big\}\Big),\quad m\in\nn,
\]
where we wrote $\|\cdot\|$ for the Euclidean norm. Applying It\^o's formula to the semimartingale $Q(\cdot)$ and recalling from \cite{RW} that the functions $\,\phi_\varepsilon (\cdot)\,$ are %chosen to be 
harmonic with respect to the generator of $B(\cdot)$, one obtains for all $\varepsilon\in(0,1)$ and $m\in\nn$:
\[
\phi_\varepsilon(Q(\tau_m))-\phi_\varepsilon(Q(0))=\int_0^{\tau_m} \nabla\phi_\varepsilon(Q(s))\cdot\mathrm{d}B(s)
+\sum_{k=1}^{n-1}\sum_{J\in\mathfrak{J}_k} \int_0^{\tau_m} v_J\cdot\nabla\phi_\varepsilon(Q(s))\,\mathrm{d}\mY^{(k)}_J\,.
\]
Here $\mathfrak{J}_k$ stands for the set of all subsets of $\{1,2,\ldots,n-1\}$ with $k$ elements and $v_J$ denotes the $J$-th column of $\mR^{(|J|)}$. Finally, taking the expectation on both sides and using the bound (18) in \cite{RW} and the bound \eqref{grad_bound} above, one ends up with
\begin{eqnarray*}
\ev[\phi_\varepsilon(Q(\tau_m))]-\ev[\phi_\varepsilon(Q(0))]
\geq-(\log\varepsilon +1) \sum_{j=1}^{n-1} c_j\,\ev\Big[\int_0^{\tau_m} \mathbf{1}_{\{\|Q(s)\|<\varepsilon\beta_j\}}\,\mathrm{d}\mY^{(1)}_j(s)\Big] \\  
-\sum_{k=1}^{n-1} C_k \sum_{J\in\mathfrak{J}_k} \ev[\mY^{(k)}_J(\tau_m)],
\end{eqnarray*}
where the positive constants $c_j$, $\beta_j$, $j=1,2,\ldots,n-1$ are defined as in \cite{RW} and the constants $C_k$, $k\in\{1,2,\ldots,n-1\}$ are as in \eqref{grad_bound}. Dividing both sides of the latter inequality by $(\log\varepsilon +1)$ and taking the limit $\varepsilon\downarrow0$ gives
\[
\lim_{\varepsilon\downarrow0} \, \sum_{j=1}^{n-1} c_j\,\ev\Big[\int_0^{\tau_m} \mathbf{1}_{\{\|Q(s)\|<\varepsilon\beta_j\}}\,\mathrm{d}\mY^{(1)}_j(s)\Big]\leq0.
\]
Thus, by Fatou's Lemma and the nonnegativity of the integrand, one can conclude
\[
\int_0^{\tau_m} \mathbf{1}_{\{Q(s)=0\}}\,\mathrm{d}\mY^{(1)}_j(s)=0,\quad j=1,2,\ldots,n-1,
\]
with probability one, and therefore also
\eq
\int_0^\infty \mathbf{1}_{\{Q(s)=0\}}\,\mathrm{d}\mY^{(1)}_j(s)=0,\quad j=1,2,\ldots,n-1, 
\en
with probability one. Finally, the backward induction argument of Lemma 5 in \cite{RW} allows us to strengthen this statement to 
\eq
\forall~~ j_0\in J\subset\{1,2,\ldots,n-1\}:\quad \int_0^\infty \mathbf{1}_{\{Q_j(s)=0,\;j\in J\}}\,\mathrm{d}\mY^{(1)}_{j_0}(s)=0,
\en
almost surely. In view of \eqref{loccond1} and \eqref{loccond2}, this finishes the proof. \ep 

%%%%%%%%%%%%%%%%%%%%%%%%%%%%%%%%%%%%%%%
\subsection{Absence of triple collisions}
%%%%%%%%%%%%%%%%%%%%%%%%%%%%%%%%%%%%%%% 

This  subsection, and the one that follows,  are devoted to the construction of the solution to the system of stochastic equations \eqref{names}, when one of the following two conditions holds:
\begin{eqnarray*}
&&\mathbf{(A)}\quad (1-q_k^-)\,\sigma_k^2\, \geq \, q_k^-\,\sigma_{k+1}^2\, ,\;\; (1-q_k^+)\,\sigma_k^2 \, \geq \, q_k^+\,\sigma_{k-1}^2\, ,\quad k=2,3,\ldots,n-1\\
&&\mathbf{(B)}\quad q_k^-\,=\,q_k^+\,=\,\left( 1+\frac{\,\sigma_{k-1}^2+\sigma_{k+1}^2\,}{2\,\sigma_k^2}\right)^{-1}\,,\quad k=2,3,\ldots,n-1.
\end{eqnarray*}
As we show below, each of these two conditions prevents collisions of three or more particles. It is   not hard to see that neither of these two conditions implies the other.

When all the collision parameters are equal to 1/2, as in (\ref{half}), condition (\textbf{B})  mandates that the graph of the variances-by-rank $\, k \mapsto \sigma^2_k\,$ be linear; thus, under   condition (\textbf{B}),  Proposition \ref{no_triple_coll} below  generalizes the results in \cite{IK} on the absence of triple collisions to situations where particles feel local-time-like drag from their immediate neighbors, when they collide with each other.

\medskip

We start with   a result  ruling out triple collisions in the case $\,n=3\,$ under the condition (\ref{concave}) below; in this three-dimensional case, condition (\ref{concave}) is weaker than each of  the conditions (\textbf{A}) and (\textbf{B}).

\begin{prop}
\label{triple_3part}
Suppose that $\,n=3$, $R_1(0)-R_3(0)>0$ and the condition 
\eq
 \label{concave}
2\, \sigma_2^2\, \geq \,q_2^-  \big(\sigma_2^2+\sigma_3^2\,\big)+q_2^+\big(\sigma_1^2+\sigma_2^2\,\big)
\en
holds. Then we have
\eq
\pp\,\big(\, \exists\,t\geq0:\; R_1(t)=R_2(t)=R_3(t)\,\big)=0\,;
\en
moreover, the converse is also true. 
\end{prop} 

\smallskip
\noindent{\bf Proof:} We consider the reflected Brownian motion $\,(R_1(\cdot) -R_2 (\cdot),\,R_2 (\cdot)-R_3(\cdot))\,$ in the nonnegative quadrant,  with the reflection matrix $\,\mathcal R=\mathbf{I}_{2}-\mathcal Q\,$ of \eqref{refl_matrix} and the covariance matrix $\mathcal A$ of \eqref{cov_matrix} with $\, n=3\,$. Moreover, let $\,\mathcal O\,$ be the $2\times 2$ orthogonal matrix such that $\, \mathcal L:=\mathcal O'\mathcal A\,\mathcal O\,$ is a diagonal matrix, where the superscript $'$ denotes transposition. Then the process
\eq
\mathcal L^{-1}\,\mathcal O\,\big(R_1(\cdot)-R_2(\cdot),\, R_2(\cdot)-R_3(\cdot)\big)'
\en  
is a reflected Brownian motion in a wedge, in the sense of \textsc{Varadhan \& Williams} (1984) in \cite{VW}. Letting 
\eq
 \label{diag}
\mathcal D\,:= \,\mathrm{diag}(\mathcal{A})
\en
 be the diagonal matrix  whose diagonal entries coincide with those of $\,\mathcal A\,$, and following the computations in subsection 3.2.1 of \cite{IK}, one concludes that the normal vectors to the two sides of the wedge are given by the columns of the matrix
\eq
\overline{\mathcal N}:=\mathcal L^{1/2}\,\mathcal O\,\mathcal D^{-1/2}\,; 
\en
whereas the reflection matrix  of the new reflected Brownian motion in the sense of \cite{VW} is given by 
\[
 \overline{\mathcal Q}\,:=\, \mathcal L^{-1/2}\,\mathcal O\,\mathcal R\,\mathcal D^{1/2}-\overline{\mathcal N}.
\]  
Furthermore, as observed in the proof of Lemma 3.2 in \cite{IK}, the corner of the wedge is attainable by the new reflected Brownian motion, if and only if the sum of the two off-diagonal entries of the matrix $\overline{\mathcal N}'\,\overline{\mathcal Q}$ is nonnegative. 

Next, we note that
\eq
\overline{\mathcal N}'\,\overline{\mathcal Q}\,=\,\mathcal D^{-1/2}\,\mathcal R\,\mathcal D^{1/2}-\overline{\mathcal N}'\,\overline{\mathcal N}.
\en
Moreover, both off-diagonal entries of $\overline{\mathcal N}'\,\overline{\mathcal N}$ are given by the negative cosine of the angle between the two sides of the wedge, which can be computed to be
\[
-\,\frac{(\mathcal A^{-1}\mathfrak{e}_1)'\,\mathfrak{e}_2}{\,((\mathcal A^{-1}\mathfrak{e}_1)'\,\mathfrak{e}_1)^{1/2} \,((\mathcal A^{-1}\mathfrak{e}_2)'\,\mathfrak{e}_2)^{1/2}\,}\,\,,
\] 
where $\,\mathfrak{e}_1,\,\mathfrak{e}_2\,$ is the canonical basis of $\,\rr^2$. Putting everything together, we can compute the sum of the two off-diagonal entries of the matrix $\,\overline{\mathcal N}'\,\overline{\mathcal Q}\,$ as
\[
-\,q_2^-\,\frac{\sqrt{\sigma_2^2+\sigma_3^2\,}}{\,\sqrt{\sigma_1^2+\sigma_2^2\,}\,}
\,-\,q_2^+\,\frac{\,\sqrt{\sigma_1^2+\sigma_2^2\,}\,}{\sqrt{\sigma_2^2+\sigma_3^2\,}}
\,+\,\frac{2\,\sigma_2^2}{\,\sqrt{(\sigma_1^2+\sigma_2^2)(\sigma_2^2+\sigma_3^2)\,}\,}\,.
\]
Simplifying this expression, we obtain the condition (\ref{concave}) of the proposition. \ep

\medskip

Now, we turn our attention to the case of general $\,n\ge 3\,$. 

\begin{prop}
 \label{no_triple_coll}
Suppose that $\, n \ge 3\,$, and that condition (\textbf{A}) or condition (\textbf{B}) holds. Then no triple collisions are possible,  that is, 
\eq
 \label{NoTriplColl}
\pp \, \big(\exists\,t\geq0,\;1\leq i<j<k\leq n:\; R_i(t)=R_j(t)=R_k(t)\big)=0\,.
\en
\end{prop}

 \smallskip

The proof   relies on  an inductive argument and the following lemma.

\begin{lemma}\label{linalg}
Let $n\geq3$. Then, condition (\textbf{A}) is equivalent to the following condition: 
\eq\label{A'cond}
\forall ~~~1\leq j,\,k\leq n-1,\;j\neq k:\;\;\big[\,\mathcal A^{-1}\big(\,\mathbf{I}_{n-1}-\mathcal Q\,\big)\big]_{jk}\geq0.
\en 
\end{lemma}

\medskip

\noindent{\bf Proof of Lemma  \ref{linalg}:} We start by recalling that the entries of $\mathcal A^{-1}$ can be computed by the formula
\eq
\mathcal A^{-1}_{jk}=(-1)^{j+k}\,\frac{\,\det (\mathcal A^{j,k})\,}{\det(\mathcal A)}\,,\quad 1\leq j,\,k\leq n-1\,,
\en
where $\mathcal A^{j,k}$ is the $(n-2)\times(n-2)$ submatrix of the symmetric matrix $\,\mathcal A\,$, that one  obtains by removing the $j$-th row and the $k$-th column from $\mathcal A\,$. Next, we introduce the notation
\eq
\Pi_{k_1}^{k_2}\,:=\,\Big(\prod_{k=k_1}^{k_2} \sigma_k^2\Big)\,\Big(\sum_{k=k_1}^{k_2} \frac{1}{\sigma_k^2}\Big)\,,\quad 1\leq k_1\leq k_2\leq n\,,
\en
and claim that the determinants $\,\det (\mathcal A^{j,k})\,$, $\,1\leq j\leq k\leq n-1\,$ are given by %the expression 
\eq
(-1)^{j+k}\,\det (\mathcal A^{j,k})\,=\, \sigma^2_1\sigma_2^2\cdots\sigma^2_n \,\Big( \sum_{\ell =1}^j \frac{1}{\sigma^2_\ell}\Big) \Big( \sum_{\ell =k+1}^n \frac{1}{\sigma^2_\ell} \Big)\,=\,\Pi_1^j\,\Big(\prod_{\ell=j+1}^k \sigma_\ell^2\Big)\,\Pi_{k+1}^n\, .
\en
This claim can be verified easily,  using induction over $n$ and distinguishing the cases $k\geq j+2$, $k=j+1$ and $k=j$. Therefore, for any fixed $1\leq j<k\leq n-1$, the inequality
\eq
 \label{A}
\big[\mathcal A^{-1}\big(\mathbf{I}_{n-1}-\mathcal Q\big)\big]_{jk}\geq0
\en      
of condition \eqref{A'cond} can be rewritten in the equivalent form
\eq
\label{Acond1}
\Pi_1^j\,\Big(\prod_{\ell=j+1}^k \sigma_\ell^2\Big)\,\Pi_{k+1}^n \,
\geq \, q_k^-\,\Pi_1^j\,\Big(\prod_{\ell=j+1}^{k-1} \sigma_\ell^2\Big)\,\Pi_k^n\,+\,(1-q_k^-)\,\Pi_1^j\,\Big(\prod_{\ell=j+1}^{k+1} \sigma_\ell^2\Big)\,\Pi_{k+2}^n\,.
\en
In view of the strict positivity of the variances, the inequality \eqref{Acond1} simplifies to
\eq
\label{Acond2}
\sigma_k^2\,\Pi_{k+1}^n\,\geq \,q_k^-\,\Pi_k^n +(1-q_k^-)\,\sigma_k^2\,\sigma_{k+1}^2\,\Pi^n_{k+2}\,. 
\en
Now, we note the following relations among $\Pi_k^n$, $\Pi_{k+1}^n\,$ and $\,\Pi_{k+2}^n\,$:
\begin{eqnarray}
%+\sigma_k^2\,\sigma_{k+2}^2\,\sigma_{k+3}^2\,\cdots\,\sigma_n^2,
&&\Pi_{k+1}^n\,=\,\Pi_{k+2}^n\,\sigma_{k+1}^2+\sigma_{k+2}^2\,\sigma_{k+3}^2\,\cdots\,\sigma_n^2\,,
\\
&&\Pi_k^n\,=\,\Pi_{k+2}^n\,\sigma_k^2\,\sigma_{k+1}^2+ \big( \sigma_{k+1}^2 + \sigma_{k}^2 \big) \,\sigma_{k+2}^2\sigma_{k+3}^2\,\cdots\,\sigma_n^2\,.
\end{eqnarray}
Plugging these into \eqref{Acond2} and simplifying further, we end up with
\eq
\label{Asimplified}
(1-q_k^-)\,\sigma_k^2 \, \geq\, q_k^-\,\sigma_{k+1}^2\,. 
\en

An analogous computation, now with $1\leq k<j\leq n-1$, shows that the inequality \eqref{A} in this case is equivalent to
\eq
\label{Asimplified2}
(1-q_j^+)\,\sigma_j^2 \, \geq\, q_j^+\,\sigma_{j-1}^2\,. 
\en
All in all, we conclude that condition \eqref{A'cond} is equivalent to condition (\textbf{A}). \ep 

\bigskip

\noindent{\bf Proof of Proposition \ref{no_triple_coll}:} We note first that, by the same application of \textsc{Girsanov}'s Theorem as in subsection 2.2 of \cite{IK}, or as in the proof of Lemmata 6 and 7 in \cite{IKS}, we need only   consider the case  $\, b_1 = b_2= \ldots = b_n =0\,$ in \eqref{ranks}. 

\medskip
\noindent $\bullet~$ 
We start by assuming that condition (\textbf{A}) is satisfied, and proceed by induction over $\,n \ge 3\,$. For $n=3$, very simple computation shows that condition (\textbf{A}) implies condition \eqref{concave}. The statement of the proposition for $\, n=3\,$ follows then directly from Proposition \ref{triple_3part}. 

\medskip

Next, we assume that $\,n\geq4\,$, and that the statement of the proposition holds under condition (\textbf{A}) for all $\,\nu=3,4,\ldots,n-1\,$ (the ``induction hypothesis"). Introducing the stopping times
\eq
 \label{tau1}
\tau_\delta\,=\,\inf\big\{t\geq0:\;\max_{k=1,2,\ldots,n-1} \big( R_k(t)- R_{k+1}(t)  \big)\leq\delta \, \big\}\,,\quad\delta>0\,,
\en 
we claim that 
\eq\label{claimdelta}
\pp\,\big(\exists\,0\leq t < \tau_\delta,\;1\leq i<j<k\leq n:\; R_i(t)=R_j(t)=R_k(t)\big)=0,\quad \delta>0. 
\en
Indeed, for any fixed $\delta>0$, the time interval $[0,\tau_\delta)$ can be written as
\eq
[0,\tau_\delta)=\bigcup_{m=0}^M \big[\tau_\delta^m,\tau_\delta^{m+1}\big)
\en
for some $M\in\nn\cup\{\infty\}$ and with stopping times $0=\tau_\delta^0<\tau_\delta^1<\ldots$ satisfying
\eq
\exists ~~~k=k(m)\in\{1,2,\ldots,n-1\}:\quad R_k(t)- R_{k+1}(t)\geq\delta\,,~\;t\in\big[\tau_\delta^m,\tau_\delta^{m+1}\big). 
\en
Moreover, on each of the time intervals $[\tau_\delta^m,\tau_\delta^{m+1})\,$, $\,0\leq m<M$, the system \eqref{ranks} splits into two subsystems of the same type which evolve independently, conditional on $R_1(\tau_\delta^m),R_2(\tau_\delta^m),\ldots,R_n(\tau_\delta^m)$. Hence, \eqref{claimdelta} is a consequence of the induction hypothesis.  

\smallskip
In view of \eqref{claimdelta}, in order to show (\ref{NoTriplColl}) and complete the induction argument, it suffices to show  
\eq
 \label{tau}
\mathbb{P}\, \Big(  \lim_{\delta\downarrow0} \tau_\delta \,= \, \infty  \Big) \,=\, 1\,.
\en
%holds with probability $1$. 
To this end, we introduce the functions  $\,F\,,\, G :\,(\rr_+)^{n-1}\backslash\{0\}\rightarrow\rr_+$ given by
\eq
F(z)\,:=\, G(z)^{(3-n)/2}\, , \qquad G(z) \,:=\,\langle \mathcal A^{-1}z,z\rangle\,. 
\en
Next, we apply the change of variable  formula of Theorem 2 in \cite{HR} to the function $F$ and the reflected Brownian motion $\, Z(\cdot)\,$ formed by the spacings (see \eqref{spacings} and \eqref{spac}), noting that the function $F$ is harmonic with respect to the generator of the  Brownian motion driving these spacings. We obtain the semimartingale decomposition 
$$
F(Z(\cdot)) \,=\, F(Z(0))+ M(\cdot) + V(\cdot)\, ,
$$
where $M(\cdot)$ is a real-valued local martingale of the form $\,M(\cdot)= \sum_{k=1}^n \int_0^{\, \cdot} \xi_k (t) \, \mathrm{d} \betab_k (t)\,$, and 
\begin{eqnarray*}
V(\cdot)\,=  
 \sum_{\ell=1}^{n-1} \int_0^{\,\cdot} \big(\partial_{z_\ell}-q_\ell^+\,\partial_{z_{\ell-1}}-q_\ell^-\,\partial_{z_{\ell+1}}\big)F(Z(s))\,\mathrm{d}\Lambda^{(\ell,\ell+1)}(t)\\
=  \frac{3-n}{2}\int_0^{\,\cdot} \big( G(Z(t))\big)^{(1-n)/2} \sum_{\ell,k=1}^{n-1}\,Z_k(t)\,
\big(\mathcal A^{-1}_{k \ell}-q_\ell^+\mathcal A^{-1}_{k\,(\ell-1)}-q_{\ell+1}^-\mathcal A^{-1}_{k\,(\ell+1)}\big)\,\mathrm{d}\Lambda^{(\ell,\ell+1)} (t)\,,  
\end{eqnarray*}
a process of finite variation on compact intervals. We also note that the processes $\, \xi_k (\, \cdot \wedge \tau_\delta)\,$ are all uniformly bounded, so the stopped local martingale $\, M (\, \cdot \wedge \tau_\delta)\,$ is in fact a martingale. 

In view of Lemma \ref{linalg}, condition (\textbf{A}) ensures that, in this last expression for the finite variation process $\,V(\cdot)\,$, the co\"efficients appearing in front of the local time processes are all non-positive, so we conclude
\eq\label{lyapunov_bound}
\forall\;t\geq0,\,\delta>0:\quad \ev \, [F(Z(t\wedge\tau_\delta))] \, \leq \, \ev \, [F(Z(0))]\,. 
\en
Moreover, if the limit $\lim_{\delta\downarrow0} \tau_\delta$ were  finite with positive probability, then for any given real number $\,c>0\,$ there would exist real numbers $\,t>0\,$ sufficiently large, and   $\,\delta>0\,$ sufficiently small, such that the value of the left-hand side in \eqref{lyapunov_bound} would exceed $\,c\,$; but this would then contradict \eqref{lyapunov_bound}. We conclude that (\ref{tau})  holds, and   thus the   proposition is established under condition (\textbf{A}).

\medskip
\noindent $\bullet~$ 
Condition (\textbf{B}) simply paraphrases the skew-symmetry condition for the $  (n-1)$-dimensional Brownian motion $\,Z (\cdot)\,$ with reflection on the faces of the nonnegative orthant in the sense of  \textsc{Harrison \& Williams} (1987) (cf.$\,$\cite{HW1}, \cite{HW2}, \cite{W} as well as section \ref{sec5} below), and the result of Proposition \ref{no_triple_coll} in this case can be found in Theorem 1.1(iii) of \cite{W}. \ep 

%%%%%%%%%%%%%%%%%%%%%%%%%%%%%%%%%%%%%%%
\subsection{Names}
%%%%%%%%%%%%%%%%%%%%%%%%%%%%%%%%%%%%%%%

We can now combine the results of the previous two sections with those  in \cite{FIK}, to give a construction of a strong solution to \eqref{names} satisfying \eqref{notsticky} and to show its pathwise uniqueness.

\begin{thm}
Suppose that condition (\textbf{A}) or condition (\textbf{B}) is satisfied. Then, the system \eqref{names} has a strong solution satisfying \eqref{notsticky}, and such a solution is pathwise unique. 
\end{thm}

\smallskip
\noindent{\bf Proof:} We start with the proof of strong existence, which proceeds by induction over $n$. As remarked in subsection \ref{sec1.1}, for $n=2$, 
the system of equations \eqref{names} is a special case of the system (4.13)--(4.14) in \cite{FIK}. Therefore, we may deduce strong existence for $n=2$ from their Theorem 4.2. Moreover, it is shown in that paper (see (6.13)-(6.15) in \cite{FIK}) that the distributions of the random variables $X_1(t)-X_2(t)$, $t\geq0$ have a density with respect to the Lebesgue measure on $\rr$ for every $\, t \in (0, \infty)\,$, and therefore
\eq
\ev\big[{\mathcal Leb}\big(\{t\geq0:\;X_1(t)=X_2(t)\}\big)\big]=0
\en
by \textsc{Fubini}'s theorem. Hence, $X_1(\cdot),\,X_2(\cdot)$ satisfy \eqref{notsticky}. 

\medskip
 
We now consider $\,n\geq 3\,$, and assume that a strong solution to the system \eqref{names} satisfying \eqref{notsticky} has already been constructed for all $\, \nu = 2,3,\ldots, n-1\,$ and all choices of drift, dispersion and collision parameters obeying condition (\textbf{A}) or (\textbf{B}) (the ``induction hypothesis''). We shall construct a strong solution $(X_1(\cdot),X_2(\cdot),\ldots,X_n(\cdot))$  of the system  \eqref{names} consecutively on the random time intervals
\eq\label{intervals}
\big[\eta^0_{2^{-p}},\eta^1_{2^{-p}}\big),\;\big[\eta^1_{2^{-p}},\eta^2_{2^{-p}}\big),\,\ldots\,,\big[\eta^{M(p)}_{2^{-p}},\eta^{M(p)}_{2^{-p}}\big),\quad p\in\nn\,,
\en
where $\eta^m_{2^{-p}}$, $m=0,1,\ldots,M(p)$, $p\in\nn$ are stopping times such that
\begin{eqnarray*}
\lim_{m\uparrow M(p)} \eta^m_{2^{-p}}\,= \,\inf \Big\{t\geq0:\;\max_{k=1,\ldots,n-1} \big( R^X_k(t) -R^X_{k+1}(t) \big) \leq 2^{-p}\Big\} = \tau_{2^{-p}}\,,\quad p\in \mathbb{N}_0\,,\\
\forall~~ p\in\nn,\,m=0,\ldots,M(p),\;~\exists~~ k=k(p,m):\;\; R^X_{k+1}(t)-R^X_k(t)\geq 2^{-p},\;~t\in[\eta^m_{2^{-p}},\eta^{m+1}_{2^{-p}}) 
\end{eqnarray*}
and $M(p)\in\nn\cup\{\infty\}$, $p\in\nn$. We have recalled here the notation of (\ref{tau1}), and the fact that  the process of ranks $(R^X_1(\cdot),R^X_2(\cdot),\ldots,R^X_n(\cdot))$ solves the system of equations \eqref{ranks} (recall the discussion in section \ref{sec1} for more details).  

\medskip

On each interval $[\eta^m_{2^{-p}},\eta^{m+1}_{2^{-p}})$, we define $(X_1(\cdot),X_2(\cdot),\ldots,X_n(\cdot))$ by letting the processes $(X_1(\cdot),X_2(\cdot),\ldots,X_{k(p,m)} (\cdot))$ evolve as a strong solution of the $k(p,m)$-dimensional system corresponding to the first $k(p,m)$ equations in \eqref{names}, started at the point $(X_1(\eta^m_{2^{-p}}),\ldots,X_{k(p,m)}(\eta^m_{2^{-p}}))\,$; and by letting $(X_{k(p,m)+1} (\cdot),\ldots,X_n(\cdot))$ evolve as a strong solution of the $(n-k(p,m))$-dimensional system corresponding to the last $(n-k(p,m))$ equations in \eqref{names}, started at $(X_{k(p,m)+1}(\eta^m_{2^{-p}}),\ldots,X_n(\eta^m_{2^{-p}}))$. Note that the strong solutions to the lower-dimensional systems exist by the induction hypothesis. The resulting process $(X_1(\cdot),X_2(\cdot),\ldots,X_n(\cdot))$ is a strong solution to the system of equations \eqref{names}, up to the random time 
\eq\label{hitting0}
\lim_{p\uparrow\infty} \lim_{m\uparrow M(p)} \eta^m_{2^{-p}} \,=\, \lim_{p\uparrow\infty}\tau_{2^{-p}} \,.
\en
We can conclude from Proposition \ref{no_triple_coll} that the quantity in \eqref{hitting0} must be equal to infinity with probability one. Thus, we have constructed a strong solution to \eqref{names} for all $t\in[0,\infty)$. Finally, it is clear that this solution satisfies \eqref{notsticky} by the induction hypothesis and its construction. 

\medskip

To prove pathwise uniqueness, we argue again by induction over $n$. For $n=2$, pathwise uniqueness is a consequence of Theorem 4.2 in \cite{FIK}. Now, let $n\geq3$ and assume that pathwise uniqueness holds for all $\, \nu = 2,3, \ldots, n-1\,$ and all choices of drift, dispersion  and collision parameters satisfying condition (\textbf{A}) or condition (\textbf{B}) (the new induction hypothesis). Next, suppose that $\,\widetilde{X} (\cdot)=\big(\widetilde{X}_1 (\cdot),\widetilde{X}_2,\ldots,\widetilde{X}_n(\cdot)\big)\,$ is another strong solution of \eqref{names} satisfying \eqref{notsticky}, defined on the same probability space as the strong solution $\,X(\cdot)=(X_1(\cdot),X_2(\cdot),\ldots,X_n (\cdot))\,$ constructed above. Considering the intervals in \eqref{intervals} consecutively and employing the induction hypothesis, we conclude that $ \, \widetilde{X} (\cdot)=X (\cdot)\,$ must hold up to the time given by \eqref{hitting0}. However, as we have seen above, the latter must be infinite with probability one, by virtue of Proposition \ref{no_triple_coll}. This yields the desired pathwise uniqueness. \ep 

%%%%%%%%%%%%%%%%%%%%%%%%%%%%%%%%%%%%%%%
\subsection{Skew-symmetry and invariant measures}
 \label{sec5}
%%%%%%%%%%%%%%%%%%%%%%%%%%%%%%%%%%%%%%% 
 
From \textsc{Ichiba et al.} (2011) (see equation (5.8) in \cite{IPBKF}), we know that the reflected Brownian motion of spacings in (\ref{spacings}) is skew-symmetric in the sense of \textsc{Harrison \& Williams} (1987) (cf.$\,$\cite{HW1}, \cite{HW2}), if the following  
condition is satisfied 
\eq
2\, \big( \,\mathcal D - \mathcal{A} \,\big) \,=\, \mathcal Q\,\mathcal D+\mathcal D\,\mathcal Q 
\en 
with the notation $\, \mathcal D=\mathrm{diag}(\mathcal{A})\,$  of (\ref{diag}). Plugging in these equations the expressions  for the matrices $\,\mathcal{A}\,$ and $\,\mathcal Q\,$ from (\ref{cov_matrix}) and (\ref{refl_matrix}), respectively, we can simplify this condition to
\eq\label{skewq}
q_k^-\,=\,q_k^+\,=\,\left( 1+\frac{\,\sigma_{k-1}^2+\sigma_{k+1}^2\,}{2\,\sigma_k^2}\right)^{-1}\,,\quad k=2,3,\ldots,n-1\,,  
\en
that is, exactly the condition ($\mathbf{ B}$). 

In view of our assumption \eqref{elastic}, the condition (\ref{skewq})  amounts to the requirement
\eq\label{skew}
\frac{2\,\sigma_2^2}{\,\sigma_1^2+\sigma_3^2\,}\,=\,\frac{\,\sigma_2^2+\sigma_4^2\,}{2\,\sigma_3^2}\,=\,\frac{2\,\sigma_4^2}{\,\sigma_3^2+\sigma_5^2\,}\,=\,\cdots\,.
\en

%\medskip

Figure 1 shows the variances $\,\sigma_k^2$, $k=1,2,\ldots,n\,$, and Figure 2 the slopes $$\, k\, \longmapsto  \, \frac{\log(\sigma^2_{k+1}-\sigma_1^2)-\log(\sigma^2_k-\sigma_1^2)}{\,\,\log k-\log (k-1)\,\,}$$ of the function $\, \log k\mapsto\log(\sigma_{k+1}^2-\sigma_1^2)\,$, for $n=100$ and a nonlinear choice of initial parameters $\sigma_1^2$, $\sigma_2^2$, $\sigma_3^2$, namely $\sigma_1^2=0.1$, $\sigma_2^2=0.11$, $\sigma_3^2=0.121$. One can see that these slopes can be made to deviate significantly from $1$ even when one only slightly perturbes a linear specification (all the slopes would be equal to $1$ in a skew-symmetric specification of a model with symmetric collisions: $q_k^{\pm}= 1/2\,$, $\,k=1,2,\ldots,n-1$). 

 \smallskip

\begin{figure}[t] 
\begin{minipage}[b]{0.45\linewidth}
\begin{center}
\includegraphics[height=2in,width=2.5in]{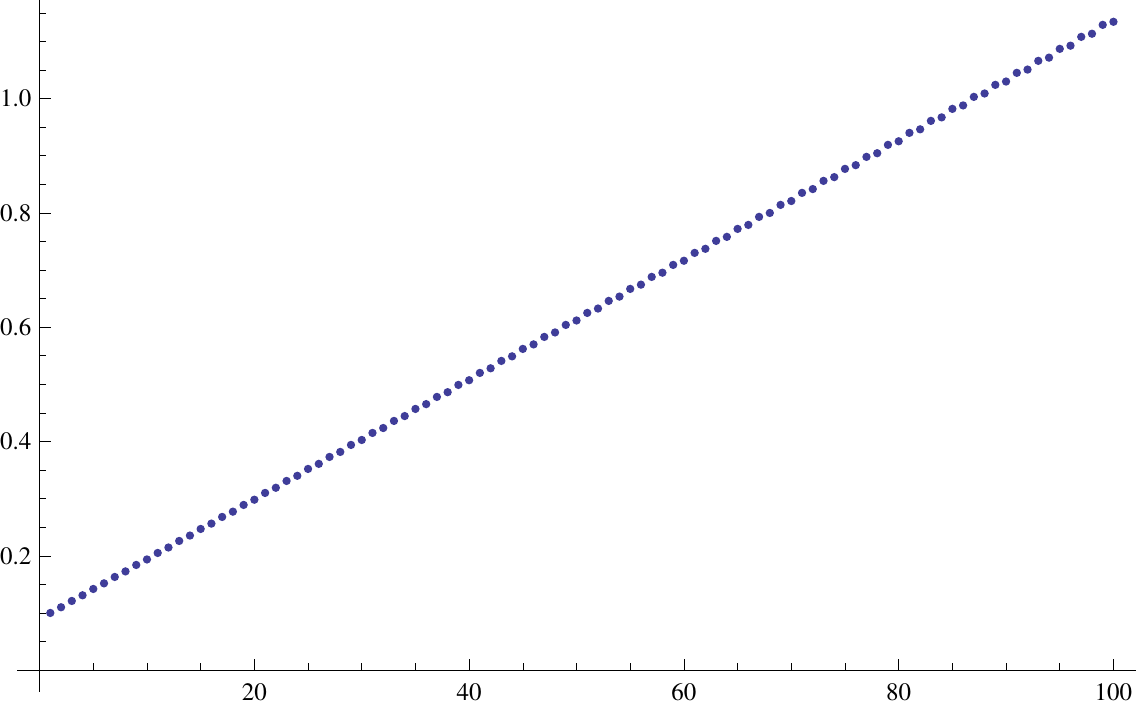}
\caption{}
\end{center}
\label{variances}
%\end{figure}
\end{minipage}
\hspace{0.5cm}
\begin{minipage}[b]{0.45\linewidth}
%\begin{figure}[t] 
\begin{center}
\includegraphics[height=2in,width=2.5in]{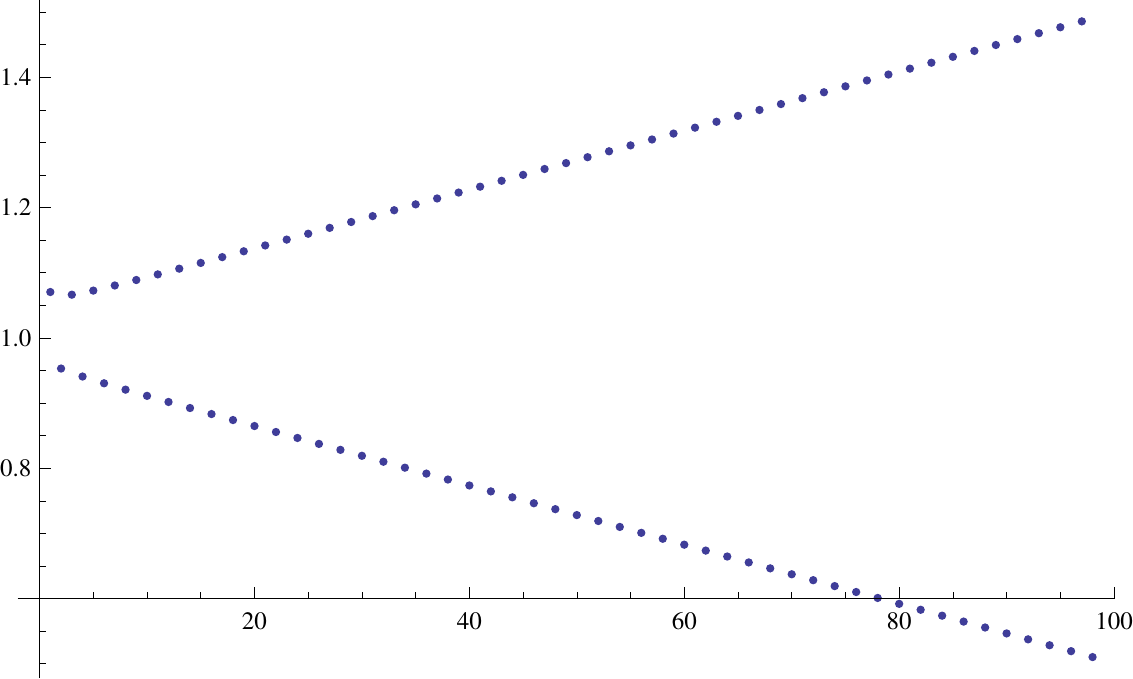}
\caption{}
\end{center}
\label{slopes}
\end{minipage}
\end{figure}

%Soumik changed dimensions+positions of figures

%\smallskip
Under the skew-symmetry condition \eqref{skew}, one can compute the invariant distribution of the spacings process \eqref{spacings}   explicitly,  using Lemma 3.6 in the dissertation of \textsc{Ichiba} (2009) \cite{I};   one ends up with a product of exponential distributions with parameter vector
\eq
\gamma\,=\,2\,\big[{\rm diag}(\mathcal{A})\big]^{-1}\, \mathcal{R}^{-1} \big(b_2-b_1,b_3-b_2,\ldots,b_n-b_{n-1}\big)'\,.
\en
Note that, by virtue of \eqref{elastic} and \eqref{skewq}, the matrix $\, \mathcal R\,$ takes the form
\eq
 \mathcal R\,=\, \mathbf{ I}_{\,n-1}-\mathcal Q\,
=\left(\begin{array}{ccccc}
1 & -q & 0 & 0 & 0 \\
-q & 1 & -(1-q) & 0 & 0  \\
0 & -(1-q) & 1 & -q & 0  \\
0 & 0 & -q & 1 & \ddots \\
0 & 0 & 0 & \ddots & \ddots
\end{array}\right),
\en
where 
\eq
q\,:=\,q_2^-\,=\,\left( 1+\frac{\,\sigma_1^2+\sigma_3^2\,}{2\,\sigma_2^2}\right)^{-1}\,.
\en

\begin{rem}
{\rm 
An interesting special case is the specification 
\eq
b_1=b_2=\ldots=b_{n-1}=0\,,\quad b_n=\mathfrak{g}\,n\;\;~~\mathrm{for\;some\;}\mathfrak{g}>0,
\en
which might be called a \textit{$q$-Atlas model} by analogy with the term Atlas model introduced by \textsc{Fernholz}  (2002) in \cite{F} and studied further by \textsc{Banner}, \textsc{Fernholz} \& \textsc{Karatzas} (2005) in \cite{BFK}. Figure 3 shows on a log-log plot the capital distribution curve
\eq
k\,\longmapsto\, \frac{e^{R_k(t)}}{\, \sum_{\ell=1}^n e^{R_\ell(t)}\,}\,,\qquad k=1,2,\ldots,100
\en
for such a $q$-Atlas model with $n=100$, $\mathfrak{g}=1$, $\sigma_1^2=0.1$, $\sigma_2^2=0.11$, $\sigma_3^2=0.121$, when the spacings process $(R_1(\cdot)-R_2(\cdot),R_2(\cdot)-R_3(\cdot),\ldots,R_{n-1}(\cdot)-R_n(\cdot))$ takes the mean value under its stationary distribution. By comparing with the plots of real-world capital disribution curves from U.S. equity market data of the Center of Research in Securities Prices (CRSP) at the University of Chicago (see Figure 5.1 on page 95 of \cite{F}), one sees that a $q$-Atlas model can capture the concave shape of the capital distribution curve, as well as its linear structure at the top. \ep}
\end{rem}

\smallskip
%Soumik reduced the dimensions 09/14
\begin{figure}[h] 
\begin{center}
\includegraphics[height=2in,width=3.5in]{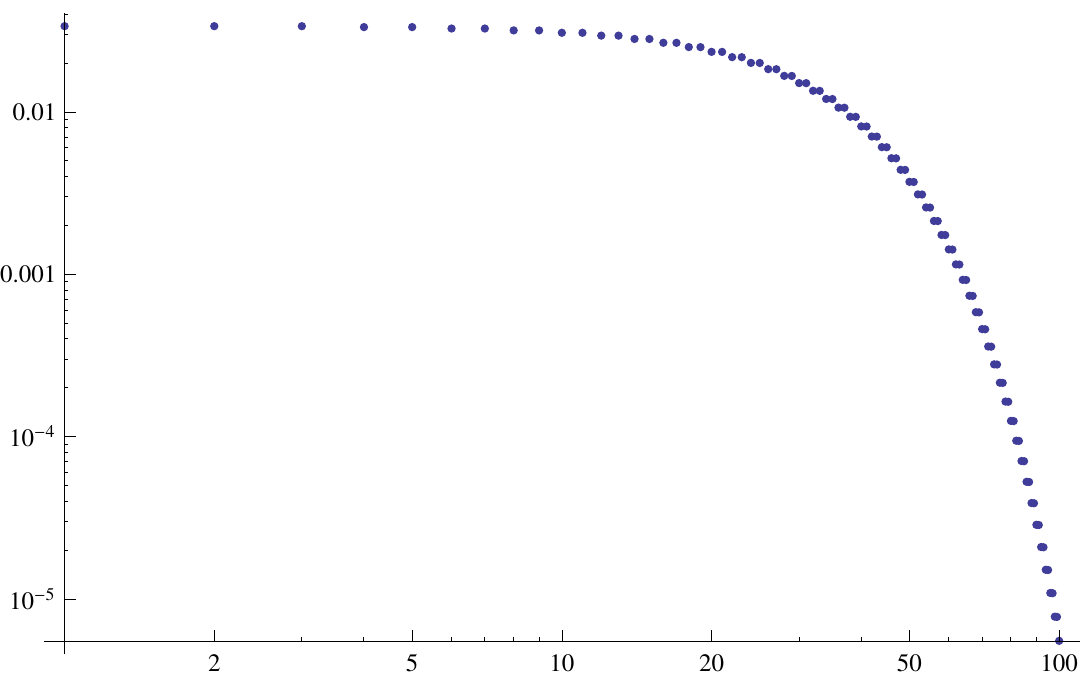}
\caption{}
\end{center}
\label{capdist}
\end{figure}

%\medskip

We also note that, since the mapping 
\eq
\big(b_2-b_1,\,b_3-b_2,\ldots,\,b_n-b_{n-1}\big)\,\longmapsto\,\gamma
\en
is bijective, one can determine the drifts up to an additive constant by fitting the vector $\gamma$ to the observed capital distribution curves, such as the ones in Figure 5.1 on page 95 of \cite{F}.

\medskip

%%%%%%%%%%%%%%Soumik-0906-Added new section%%%%%%%%%%%
\section{Scaling limit of asymmetrically colliding random walks}
%%%%%%%%%%%%%%%%%%%%%%%%%%%%%%%%%%%%%%%%%

Consider again the processes corresponding to the system of stochastic equations \eqref{ranks}. Our objective in this section is to show that such processes with asymmetric local time components arise as scaling limits of random walks with asymmetric interactions upon ``collision''. We start with the following informal description.

\medskip

Consider an $n$-dimensional continuous time jump process $\rwalk(\cdot)=(\rwalk_1(\cdot),\ldots,\rwalk_n(\cdot))$ on the wedge
\[
\wedgen_n = \left\{ (z_1, \ldots, z_n)\in \mathbb{Z}^n:\; z_1 \ge z_2 \ge \ldots \ge z_n \right\}.
\]
The process registers the positions on the integer lattice of $n$ particles that always maintain their order. We will label the particles  by the elements of the set $[n]=\{1,2,\ldots,n\}$ where $1$ refers to the rightmost particle and $n$ refers to the leftmost. 

\medskip

The movement of particles can be described as follows. Consider nonnegative parameters $\{a_k,b_k,\theta_k^L,\theta_k^R,\;k\in[n]\}$. If there are no other particles at the current site of particle $k$, then this particle moves to the right at rate $a_k$ and to the left at rate $b_k$, independently of every other particle. On the other hand, suppose that  particles $k,k+1,\ldots,k+\ell$ are currently at the same location (a phenomenon we call ``collision"), and this is the maximum length of the tie in the sense that the $(k-1)$-st and the $(k+\ell+1)$-st particles are not at that site. Then, at rate $\,\theta^R_{k}a_k\,$, particle $k$ jumps to the right;  and, at rate $\,\theta_{k+\ell}^L b_{k+\ell}\,$, particle $(k+\ell)$ jumps to the left, while the particles in between do not move. 

\medskip

We wish to take a diffusion limit of such systems. Before we go on, let us give a few examples which display the variety of behavior that we can expect from such particle systems. 

\medskip
\noindent\textbf{Example 1.} Suppose each particle evolves according to a Poisson process with parameter $a$ until a collision occurs. Thus  every $a_k \equiv a\,$, every $b_k \equiv 0\,$, and the values of $\theta^L_k$, $k\in[n]$ are irrelevant. 

Consider, first, $\,\theta^{R}_k\equiv2$, $k\in[n]\,$. This corresponds  via relabeling  to the case when the Poisson particles are moving independently and are allowed to pass each other. On the other hand, when each $\theta_k^R\equiv1$, the higher particle blocks the movement of the lower one when they collide. The resulting particle system then evolves as the well-known TASEP process (except that in TASEP collisions occur when particles are at adjacent sites instead of the same site).  

\medskip
\noindent\textbf{Example 2.} Suppose the particles move according to i.i.d. simple symmetric random walks until collision. That is, $a_k\equiv b_k\equiv1$, $k\in[n]$. Moreover, we take $\theta_k^{L}\equiv\theta_k^R \equiv \theta$, $k\in[n]$. Thus, when a particle is involved in a collision, this simply changes the rate of its next jump. The case $\theta=2$ yields the ordered system of i.i.d. random walks. When $\theta\approx0$, the particles get ``sticky'', while when $\theta$ is very large, the particles can be thought of as repelling one another. 

\medskip

\noindent\textbf{Example 3.} Consider, as before, particles moving according to i.i.d. simple symmetric random walks until collision. That is, $a_k\equiv b_k\equiv 1$, $k\in[n]$. Fix an $\ell\in [n]$, say $\ell=2$. Choose the $\theta_2^{L},\,\theta^R_2$ parameters to be much larger than the rest of the collision parameters. Thus, whenever particle $1$ or particle $3$ collides with particle $2$, the latter will almost immediately move away. Hence, particle $2$ will remain \textit{invisible} to its neighbors. 
 
 \smallskip

Now, consider the triple collisions of particles $1$, $2$ and $3$. Note that, with $\theta_2^{L}$ and $\theta^R_2$ being very large, such triple collisions will happen for about the same duration of time as collisions between particle $1$ and particle $3$ in the particle system obtained by removing particle $2$. In other words, if one removes particle $2$, this should not influence the behavior of the other particles in a significant way.   

\medskip

\noindent\textbf{Example 4.} We let the particles evolve according to i.i.d. standard Poisson processes until collision, as in Example 1. However, now we set $\theta_2^R\equiv0$, with all other collision parameters being positive. Then  particle $2$ and particle $3$ freeze forever the first time they collide. Thereafter, particle $1$ moves as a standard Poisson process independently from the rest of the particle system, while particles $4,5,\ldots,n$ eventually all coalesce at the site of particle $2$.  

\medskip
%Soumik
%Yannis
\noindent
{\it Physical heuristics:} To develop a feel for these processes, let us discuss briefly the mechanics involved in these collisions. In classical models of collision (see \cite{H}), the particles behave as hard billiard balls of infinitesimally small radius. This models \textit{elastic collision} which we now explain. 

Suppose two particles collide on the line. Particle $1$, on the right, has mass $m_1$ and velocity $u_1$ right before the collision. Particle $2$, coming from the left, has mass $m_2$ and velocity $u_2$. In elastic collisions the momentum and the kinetic energy are preserved before and after the collision. These two preserved quantities uniquely determine the velocities $v_1$ and $v_2$ of the two particles right after the collision, i.e.,
\[
v_1=\frac{u_1(m_1-m_2) + 2m_2 u_2}{m_1 + m_2}, \quad v_2= \frac{u_2(m_2-m_1) + 2m_1 u_1}{m_1 + m_2}.  
\] 

Consider again our jump processes. We follow the heuristics of \cite{H} and assume that particles jump together at discrete times. A similar but slightly more lengthy analysis can be done for continuous time,  when we let one particle stay still and get hit by the other particle. We consider the Einsteinian viewpoint, that the $n$ particles are bombarded on all sides by other small particles which lead to their random motions. Suppose particles $k$ and $k+1$ are adjacent and are of masses $m_1$ and $m_2$ respectively. We will think of the rate of jumps as the speed in the appropriate direction. 

In the next small time interval, these two particles either   do not collide, and move away from one another at velocity $a_k$ and $-b_{k+1}$ respectively; or they collide, and $\,u_1=-b_{k}\,$, $\,u_2=a_{k+1}\,$. Thus, if the collision is elastic, we will observe that the velocity (total jump rate to the right)  for particle $k$ is
\[
\delta_k:=a_k + \frac{-b_k(m_1-m_2) + 2 m_2 a_{k+1} }{m_1 + m_2},
\]
and the total jump rate to the left for particle $(k+1)$ is
\[
\delta_{k+1}:=-b_{k+1} + \frac{a_{k+1}(m_2-m_1)- 2m_1 b_k }{m_1 + m_2}.
\]
Clearly, unless we impose specific constraints on the parameters, these are not equal to $\,\theta^R_k a_k\,$ and $\,-\theta^L_{k+1}b_{k+1}\,$, respectively. Hence,  these collisions are not elastic in general.  
%Soumik 09/29

Certain special cases are worth mentioning. Suppose all masses are equal to one. Then, under elastic collision, $\delta_k=a_k+ a_{k+1}$ and $\delta_{k+1}=-b_k-b_{k+1}$. This is the case when particles exchange velocities and can be thought of, via relabeling, as crossing over. A specific choice of $\theta^R_k$ and $\theta^L_{k+1}$ would capture this scenario as in the case of $\theta^R_k=2$ in Example 1 and $\theta=2$ in Example 2 above.  

On the other hand, suppose that all $\,a_k$ and $\, b_k$ are equal to 1. Then, under elastic collision, $\theta^R_k=\delta_k= 4m_2/(m_1 + m_2)$ and $\theta^L_{k+1}=\delta_{k+1}=4m_1/(m_1 + m_2)$. In particular, if one of $\theta^R_k + \theta^L_{k+1}$ is greater (or, less) than $4$, the mechanics generates excess energy (or absorbs energy) that cannot be explained by elastic collision. 

\smallskip
Under suitable assumptions on the parameters,  we establish in Theorem \ref{thm:conv_full} diffusion limits similar to \eqref{ranks} for these particle systems. There is an apparent paradox here. Consider that diffusion limit for $a=b=\sigma_k^R=\sigma_k^L=1$ for all $k\,$; somewhat surprisingly, the limit turns out to depend  only on the ratio $\,\theta^R_k/\theta^L_{k+1}\,$. In other words, the diffusion limit is the same whether $\theta^R_k + \theta^L_{k+1}$ is equal to $4$ or not. {\it The collisions among the limiting diffusion particles are always elastic,} and the quantities $q_k^-$'s and $q_k^+$'s can be thought of as the proportions of total mass shared by the colliding particles. This is   a consequence of the fact that the occupation time of collisions has Lebesgue measure zero in the limit, as will be made clear in the proof. 

To get a true inelastic limit, one has to let $\,\theta^R_k, \,\theta^L_k\,$ go to zero suitably in the diffusion scaling; then, presumably, one would obtain \textit{sticky} colliding Brownian particles.  %We plan to investigate this situation in the future.

\subsection{The modified Skorokhod problem} 

We start with a set of parameters: $a,b, (\lambda^L_k, \sigma^L_k, \theta^L_k,\; k\in [n])$, $(\lambda^R_k, \sigma^R_k, \theta^R_k,\; k\in [n])$. At this point we only assume that $\sigma^L_k$, $\sigma^R_k$ are stricly positive and $\theta^L_k$, $\theta^R_k$ are nonnegative for every $k$.  

\smallskip

Taking a diffusion limit requires considering a sequence of interacting jump processes as described above. We shall generalize the setup by allowing non-exponential waiting times for the jumps. Let the sequence of interacting jump processes be indexed by $N>0$. We fix a value of $N$ and a probability space rich enough to support mutually independent sequences of i.i.d. random variables $(u_k^L(i),\;i\in\NN)$, $k\in[n]$ and $(u_k^R(i),\;i\in\NN)$, $k\in[n]$, all taking only positive values. These random variables denote the inter-jump times of the particles (the superscripts ``$L$'' and ``$R$'' standing for leftward and rightward jumps, respectively). Assume that, for any fixed $k\in[n]$, 
\eq\label{eq:meanvar}
\begin{split}
\ev\left[u_k^L(1)\right] &= \left(b+ \frac{\lambda^L_k}{\sqrt{N}}\right)^{-1},\quad \mathrm{Var}\left[u_k^L(1)\right]=\left(\sigma^L_k\right)^2,\\
\ev\left[u_k^R(1)\right] &= \left(a+ \frac{\lambda^R_k}{\sqrt{N}}\right)^{-1},\quad \mathrm{Var}\left[u_k^R(1)\right]=\left(\sigma^R_k\right)^2.
\end{split}
\en

Next, we define the corresponding partial sum processes
\[
\begin{split}
U_k^L(0)&=0, \quad U_k^L(j)=\sum_{i=1}^j u_k^L(i), \quad j\in\NN,\quad k\in[n],\\ 
U_k^R(0)&=0,\quad U_k^R(j)= \sum_{i=1}^j u_k^R(i), \quad j\in\NN,\quad k\in[n],
\end{split}
\]
and the corresponding renewal processes
\[
\begin{split}
S_k^L(t)&=\max\left\{j\ge 0:\;U_k^L(j) \le t \right\}, \quad t\geq0,\quad k\in[n],\\
S_k^R(t)&=\max\left\{j\ge 0:\;U_k^R(j) \le t \right\}, \quad t\geq0,\quad k\in[n].
\end{split}
\]
Finally, we denote by $(\mathcal{F}_t\,, ~t\ge0)$ the filtration generated by the processes $S^L_k (\cdot)$, $S^R_k (\cdot)$, $k\in[n]$. Informally, for each $k\in[n]$, the process $S_k^L(\cdot)$ (resp., $S_k^R(\cdot)$) records the leftward (respectively, rightward) jumps of the $k$-th particle from the right, as long as this  particle  is not involved in a collision. 

\smallskip
To describe the effect of collisions we shall use a stochastic time change. The following lemma encapsulates the idea that the leftward (respectively, rightward) movement   for the  of particle $k$ can either be blocked, or proceed at a different rate, depending on whether there is a collision with particle $(k+1)$ (resp., particle $(k-1)$). Such results are standard in the Queueing Theory literature (e.g., section 2 in the seminal article \cite{R} by \textsc{Reiman} (1984)), so we omit the proof.  

\begin{lemma}\label{lem:existence}
For every $N\in \mathbb{N}$, and any $(\gamma_1,\gamma_2,\ldots,\gamma_n)\in\wedgen_n\,$, there exists a system of jump processes $\Gamma(\cdot)\equiv(\Gamma_k(\cdot),\;k\in [n])$ taking values in $\wedgen_n$ and progressively measurable with respect to $\mcal{F}_t$, $t\ge 0$, that satisfies the following set of equations pathwise:
\eq\label{eq:existence} 
\begin{split}
\Gamma_k(t)&=\gamma_k-S_k^L\left(T^L_k(t)\right)+S_k^R\left(T^R_k(t)\right),\quad \text{where}\\
T^L_k(t)&=\int_0^t \mathbf{1}_{\{Q_{k-1}(s)>0\}}\,\mathrm{d}T_k(s)+\theta^L_k \int_0^t \mathbf{1}_{\{Q_{k-1}(s)=0\}}\,\mathrm{d}T_k(s), \quad \text{and}\\
T^R_k(t)&=\int_0^t \mathbf{1}_{\{Q_k(s)>0\}}\,\mathrm{d}T_{k-1}(s)+\theta^R_k \int_0^t \mathbf{1}_{\{Q_k(s)=0\}}\,\mathrm{d}T_{k-1}(s),
\end{split}
\en
$k\in[n]$. We have denoted here  by 
\eq
Q(\cdot) \,\equiv \, \big(Q_k(\cdot),\;k\in[n-1])\, \equiv \, \big(\Gamma_k(\cdot)-\Gamma_{k+1}(\cdot),\;k\in[n-1]\big) 
\en
the process of gaps, and   have set
\[
T_k(t)= \int_0^t \mathbf{1}_{\{Q_k(s)>0\}}\,\mathrm{d}s, \quad k=0,1,\ldots,n,   
\]
with the convention $\mathbf{1}_{\{Q_0(\cdot)>0\}}\equiv\mathbf{1}_{\{Q_n(\cdot)>0\}}\equiv1$. 
\end{lemma}

\smallskip

It should be noted that, although we have suppressed the index $N$ from the notation in the lemma above, this parameter determines the drifts of the 
%Soumik
co\"ordinate processes.
The key to passing to the scaling limit is to understand the time-changes involved. Our strategy is the following: (i) express the process of gaps as a Skorokhod map applied to a suitable ``noise process''; (ii) show that the sequence of distributions of gaps is tight; (iii) and finally, show that   tightness implies the convergence of the appropriately rescaled process $\Gamma(\cdot)$ to a semimartingale of the type described in \eqref{ranks}.   

\medskip

We start by analyzing the process of gaps. To this end, we define the centered processes 
\[
\centS^L_k(t)=S_k^L(t)-bt,\;t\geq0,\quad \centS^R_k(t)=S_k^R(t)-at,\;t\geq0
\]
for all $k\in[n]$. Also, we define the following processes measuring the time spent in the various collisions:
\[
\begin{split}
&I_k(t)=t-T_k(t),\;t\geq0,\quad k=0,1,\ldots,n,\\
&I_{k,k+1}(t)=\int_0^t \mathbf{1}_{\{Q_k(s)=Q_{k+1}(s)=0\}}\,\mathrm{d}s,\quad k=0,1,\ldots,n-1.
\end{split}
\]
With this  notation,
\begin{eqnarray*}
Q_k(t)&=&\gamma_k-\gamma_{k+1}-S_k^L\left(T^L_k(t)\right)+S_k^R\left(T^R_k(t)\right)+S_{k+1}^L\left(T^L_{k+1}(t)\right)
-S_{k+1}^R\left(T^R_{k+1}(t)\right) \\
&=&\gamma_k-\gamma_{k+1}-\centS_k^L\left(T^L_k(t)\right)+\centS_k^R\left(T^R_k(t)\right)+\centS_{k+1}^L\left(T^L_{k+1}(t)\right)-\centS_{k+1}^R\left(T^R_{k+1}(t)\right) \\
&&-b\left(T_k^L(t)-T_{k+1}^L(t)\right)+a\left(T_k^R(t)-T_{k+1}^R(t)\right),\quad k\in[n-1].
\end{eqnarray*}
Moreover,
\eq\label{eq:whatisb}
\begin{split}
T_k^L(t)&=T_k(t)+\left(\theta_k^L-1\right)\int_0^t \mathbf{1}_{\{Q_{k-1}(s)=0,\,Q_k(s)>0\}}\,\mathrm{d}s \\
&=t-I_k(t)+\left(\theta_k^L-1\right)\left(I_{k-1}(t)-\int_0^t \mathbf{1}_{\{Q_{k-1}(s)=Q_k(s)=0\}}\,\mathrm{d}s\right) \\
&=t-I_k(t)+\left(\theta_k^L-1\right)I_{k-1}(t)-\left(\theta_k^L-1\right)I_{k-1,k}(t),\quad t\geq0 
\end{split}
\en
and 
\eq\label{eq:whatisb2}
T_k^R(t)=t-I_{k-1}(t)+\left(\theta^R_k-1\right)I_k(t)-\left(\theta^R_k-1\right)I_{k-1,k}(t),\quad t\geq0
\en
for all $k\in[n]$. 

Next, we introduce the   process $\centX(\cdot)\equiv(\centX_k(\cdot),\;k\in[n-1])$ via
\[
\centX_k(\cdot)=-\centS_k^L\left(T^L_k(\cdot)\right)+\centS_k^R\left(T^R_k(\cdot)\right)+\centS_{k+1}^L\left(T^L_{k+1}(\cdot)\right)-  \centS_{k+1}^R\left(T^R_{k+1}(\cdot)\right),
\]
$k\in[n]$; this   will play the r\^ole of   ``noise process'', to which a Skorokhod map will be applied to obtain the process of gaps $Q(\cdot)\equiv(Q_k(\cdot),\;k\in[n-1])$. 

Combining everything so far, we obtain the representation for  these gaps
\eq\label{eq:expressqueue}
\begin{split}
Q_k(t)=&\,\gamma_k-\gamma_{k+1}+\centX_k(t)\\
&-b\left(I_{k+1}(t)-I_k(t)\right)-b\left(\theta^L_k-1\right)I_{k-1}(t)+b\left(\theta^L_k-1\right)I_{k-1,k}(t)\\
&+b\left(\theta^L_{k+1}-1\right)I_k(t)-b\left(\theta^L_{k+1}-1\right)I_{k,k+1}(t)\\
&+a\left(I_k(t)-I_{k-1}(t)\right)+a\left(\theta^R_k-1\right)I_k(t)-a\left(\theta^R_k-1\right)I_{k-1,k}(t)\\
&-a\left(\theta^R_{k+1}-1\right)I_{k+1}(t)+a\left(\theta^R_{k+1}-1\right) I_{k,k+1}(t),\quad k\in[n].
\end{split}
\en

\smallskip

Now, it is very convenient to write this equation in a more transparent, matrix-vector notation, so we introduce suitable ``reflection matrices''. Let $\mR$ be the $(n-1)\times (n-1)$ matrix such that, for every $k\in[n-1]$, all entries in the $k$-th row of $\mR$ are zero, except the 
$(k-1)$-st, the $k$-th and the $(k+1)$-st. These are given by
\eq\label{eq:whatisr}
\mR_{k,k-1}=-\frac{a+b\left(\theta^L_k-1\right)}{a\theta^R_{k-1}+b\theta^L_k}, \quad 
\mR_{k,k}=1, \quad 
\mR_{k,k+1}=-\frac{b+a\left(\theta^R_{k+1} - 1\right)}{a\theta^R_{k+1}+b\theta^L_{k+2}}.
\en
Note that the columns of the matrix $\, \mathbf{I}_{n-1}-\mR\,$ add up to one. Similarly, define the $(n-1)\times(n-2)$ matrix $\modR$ such that, for each $k\in[n-1]$, all entries in the $k$-th row of $\modR$ are zero, except for the $(k-1)$-st and the $k$-th, which are given as
\eq\label{eq:whatisrtilde}
\modR_{k,k-1}=-a\left(\theta^R_k-1\right)+b\left(\theta^L_k-1\right), \quad 
\modR_{k,k}=a\left(\theta^R_{k+1}-1\right)-b\left(\theta^L_{k+1}-1\right). 
\en
Finally, we introduce the stochastic processes
\eq\label{eq:whatisy}
\mY_k(\cdot)\equiv\left(a\theta^R_k+b\theta^L_{k+1}\right) I_k(\cdot),\; k\in [n-1], \quad \modY_k(\cdot)\equiv I_{k,k+1}(\cdot),\; k\in [n-2]. 
\en
With this  notation, we have the following cleaner matrix-vector analogue of \eqref{eq:expressqueue}:
\eq\label{eq:expressq2}
Q(\cdot)=Q(0)+\centX(\cdot)+\mR\,\mY(\cdot)+\modR\,\modY(\cdot).
\en
We shall refer to this equation as the {\it modified Skorokhod representation}. 

%%%%%%%%%%%%%%%%%%%%%%%%%%%%%%%%%%%%
\subsection{The diffusion limit}
%%%%%%%%%%%%%%%%%%%%%%%%%%%%%%%%%%%%

To be able to pass to the diffusion limit, we make the following assumptions on the parameters. 

\begin{asmp}\label{coeff_asmp}
{\it Assume that $a>0$, $b>0$,
\begin{enumerate}
\item[(i)] the entries of $\mR$ satisfy 
\[
\frac{a+b\left(\theta^L_k-1\right)}{a\theta^R_{k-1}+b\theta^L_k}\in(0,\infty)\,,\;\;\;~~
\frac{b+a\left(\theta^R_{k+1}-1\right)}{a\theta^R_{k+1}+b\theta^L_{k+2}}\in(0,\infty)\,,\quad  \quad k\in[n-1]\,\mathrm{;}
\]
\item[(ii)] there is a $k_0\in[n+1]$ such that the numbers $\,a\left(\theta_k^R-1\right)-b\left(\theta_k^L -1\right)$, $k\in[n]$ are nonpositive for all $k<k_0$ and nonnegative for all $k\geq k_0$;
\item[(iii)] for some $\varepsilon>0$, we have  
\[
\sup_{N\in \NN}\;\max_{k\in [n]}\;\left(\ev\left[u_k^L(1)^{2+\varepsilon}\right]+\ev\left[u_k^R(1)^{2+\varepsilon}\right]\right) < \infty.
\] 
\end{enumerate}
}
\end{asmp}

\smallskip

A simple case, in which parts (i) and (ii) of Assumption \ref{coeff_asmp} are satisfied, is given by $\theta_k^L=\theta_k^R\geq1$, $k\in[n]$. Next, for each $m\in\NN$, we write $D^m[0,\infty)$ for the space of right-continuous $\rr^m$-valued functions on $[0,\infty)$ having left limits, endowed with the topology of uniform convergence on compact sets. 

On the strength of part (iii) of the above assumption, the following result follows from Theorem 14.6 in \cite{B99} (note that one can improve the topology in the conclusion of Lemma \ref{lem:BMconv} to the locally uniform topology, by noting the path continuity of the limit process, modifying the paths of the jump processes to continuous, piecewise linear functions, and using the fact that the Skorokhod topology relativized to the space of continuous functions coincides with the locally uniform topology there).  

\begin{lemma}
 \label{lem:BMconv}
The distribution of the process
\[
\left(\frac{1}{\sqrt{N}}\,\centS_k^L\left(Nt\right),\,t\geq0,\;\;\frac{1}{\sqrt{N}}\,\centS_k^R\left(Nt\right),\,t\geq 0,\quad k\in[n]\right)
\]
converges weakly in $D^{2n}[0,\infty)$ to the law of the vector of independent processes $\left(Z_k^L(\cdot),\;Z_k^R(\cdot),\;k\in[n]\right)$. Here, for each $k\in[n]$, $Z_k^L(\cdot)$ is a Brownian motion with drift co\"efficient $\,\lambda_k^L\,$ and diffusion co\"efficient $\,b^{3/2}\sigma_k^L\,$, while $Z_k^R(\cdot)$ is a Brownian motion with drift co\"efficient $\,\lambda_k^R\,$ and diffusion co\"efficient $\,a^{3/2}\sigma_k^R\,$. 
\end{lemma}

Next, we introduce the rescaled versions of the gap process by
\[
Q^N(\cdot)\equiv\left(Q^N_k(\cdot),\;k\in[n-1]\right)\equiv\left(\frac{1}{\sqrt{N}}\,Q_k\left(Nt\right),\;t\geq0,\quad k\in [n-1]\right)
\]
and define $\mY^N(\cdot)$, $\modY^N(\cdot)$ accordingly. Then, 
\eq\label{eq:expressq3}
Q^N(\cdot)=Q^N(0)+\centX^N(\cdot)+\mR\,\mY^N(\cdot)+\modR\,\modY^N(\cdot),
\en
where
\eq\label{eq:whatisxN}
\begin{split}
\centX^N_k(t)\equiv\frac{-1}{\sqrt{N}}\left[\centS_k^L\left(T^L_k(Nt)\right)-\centS_k^R\left(T^R_k(Nt)\right)-\centS_{k+1}^L\left(T^L_{k+1}(Nt)\right)+\centS_{k+1}^R\left(T^R_{k+1}(Nt)\right)\right],
\end{split}
\en
$t\geq0$, $k\in[n-1]$. We can   state now our main limit theorem for the process of gaps. 
\begin{thm}\label{thm_conv}
Suppose Assumption \ref{coeff_asmp} holds, and that $\lim_{N\rightarrow\infty} Q^N(0)=\xi(0)$ in distribution for some random vector $\xi(0)$. Let $\xi(\cdot)$ be a reflected Brownian motion in $(\rr_+)^{n-1}$ with initial condition $\xi(0)$,   drift vector 
\eq
\mathfrak{b}=\left(-\lambda^L_k+\lambda^R_k+\lambda^L_{k+1}-\lambda^R_{k+1},\;k\in[n-1]\right),
\en
 diffusion matrix $\,\mA=(\mA_{k,\ell})_{1\leq k,\ell\leq n-1}\,$ given by 
\eq
\mA_{k,l}=
\begin{cases}
-a^3\left(\sigma_{k+1}^R\right)^2-b^3\left(\sigma_{k+1}^L\right)^2 & \mathrm{if}\;\;\; \ell=k+1,\; k\in [n-2],\\
-a^3\left(\sigma_k^R\right)^2-b^3\left(\sigma_k^L\right)^2 & \mathrm{if}\;\;\; \ell=k-1,\; \ell\in [n-2],\\
a^3\left(\sigma_k^R\right)^2+a^3\left(\sigma_{k+1}^R\right)^2 +b^3\left(\sigma_k^L\right)^2+b^3\left(\sigma_{k+1}^L\right)^2 
& \mathrm{if}\;\;\;k=\ell\in [n-1],\\
0 & \mathrm{otherwise,}
\end{cases}
\en
and reflection matrix $\mR$ as in \eqref{eq:whatisr}. 

Then, the processes $Q^N(\cdot)$, $N>0$ converge in distribution to $\xi(\cdot)$, as $N\rightarrow\infty$, in $D^{n-1}[0,\infty)$. Moreover, the processes
\eq\label{eq:localtime}
\left(\mathfrak{Y}^N_k(\cdot),\;k\in[n-1]\right)
=\left((a\theta_k^R+b\theta_k^L)\int_0^\cdot \mathbf{1}_{\{Q^N_k(s)=0\}}\,\mathrm{d}s,\;k\in[n-1]\right)
\en  
converge in distribution in $D^{n-1}[0,\infty)$ to the process of local times accumulated by $\xi(\cdot)$ on the respective faces of the orthant $(\rr_+)^{n-1}$, and the processes 
\eq\label{eq:localtime2}
\left(\widetilde{\mathfrak{Y}}^N_k(\cdot),\;k\in[n-2]\right)
=\left(\int_0^\cdot \mathbf{1}_{\{Q^N_k(s)=Q^N_{k+1}(s)=0\}}\,\mathrm{d}s,\;k\in[n-2]\right)
\en  
tend to zero in distribution in $D^{n-2}[0,\infty)$. 
\end{thm}

We remark at this point that the limit process can be viewed as the process of gaps for a semimartingale as in \eqref{ranks} upon the appropriate identification of parameters (see Theorem \ref{thm:conv_full} below for the details). The proof of Theorem \ref{thm_conv} relies heavily on an extension of the invariance principle of \textsc{Williams} \cite{W1}, which we introduce and prove in the following subsection.  

%%%%%%%%%%%%%%%%%%%%%%%%%%%%%%%%
\subsection{Oscillation estimates for modified Skorokhod maps} \label{subs_Gip}
%%%%%%%%%%%%%%%%%%%%%%%%%%%%%%%%

In this subsection, we consider families $(Q^N(\cdot),\centX^N(\cdot),\mY^N(\cdot),\modY^N(\cdot))$, $N>0$ of processes with right-continuous paths having left limits, which satisfy the following properties. The processes $Q^N(\cdot)$, $N>0$ take values in the orthant $(\rr_+)^{n-1}$; the components of the processes $\mY^N(\cdot)$, $N>0$ and $\modY^N(\cdot)$, $N>0$ are non-decreasing and start at $0$; the modified Skorokhod representation
\eq\label{skor_eq}
Q^N(\cdot)=\centX^N(\cdot)+\mR\,\mY^N(\cdot)+\modR\,\modY^N(\cdot)
\en
holds for each $N>0$, where $\mR$ and $\modR$ are as in \eqref{eq:whatisr} and \eqref{eq:whatisrtilde}, respectively; and the following are true:
\begin{eqnarray}
&&\;\;\int_0^\infty \mathbf{1}_{\{Q^N_k(s)>0\}}\,\mathrm{d}\mY^N_k(s)=0,\quad k\in[n-1], \\
&&\;\;\int_0^\infty \big(\mathbf{1}_{\{Q^N_k(s)>0\}}+\mathbf{1}_{\{Q^N_{k+1}(s)>0\}}\big)\,\mathrm{d}\modY^N_k(s)=0, \quad k\in[n-2], \\
&&\;\;\forall\,0\leq s<t:\;\modY^N_k(t)-\modY^N_k(s)\leq\min\left(\mY^N_k(t)-\mY^N_k(s),\mY^N_{k+1}(t)-\mY^N_{k+1}(s)\right).
\end{eqnarray}
In this situation, we have the following extension of the invariance principle of \cite{W1}. 

\begin{prop}\label{inv_prin}
Under parts (i) and (ii) of Assumption \ref{coeff_asmp}, suppose that the initial conditions $Q^N(0)=X^N(0)$, $N>0$ converge in distribution to a random vector $\xi(0)$. Then:
\begin{enumerate}[(a)]
\item If the family of processes $\centX^N(\cdot)$, $N>0$ is tight on $D^{n-1}[0,\infty)$, then the same is true for $\big(Q^N(\cdot),\centX^N(\cdot),\mY^N(\cdot),\modY^N(\cdot)\big)$, $N>0$ on $D^{4n-5}[0,\infty)$.
\item In the situation of part (a), any limit point $\big(Q^\infty (\cdot),\centX^\infty (\cdot),\mY^\infty (\cdot),\modY^\infty (\cdot) \big)$ with continuous paths satisfies
	\begin{enumerate}[(i)]
	\item $Q^\infty(t)=\centX^\infty(t)+\mR\,\mY^\infty(t)+\modR\,\modY^\infty(t)\in(\rr_+)^{n-1}$, $t\geq0$. 
	\item The components of $\,\mY^\infty (\cdot)\,$ and $\,\modY^\infty (\cdot)\,$ are non-decreasing, and satisfy $\,\mY^\infty(0)=0$, $\,\modY^\infty(0)=0\,$ with probability $1$. 
	\item The following two identities hold with probability  $1 :$ 
	\begin{eqnarray}
	&&\int_0^\infty \mathbf{1}_{\{Q^\infty_k(s)>0\}}\,\mathrm{d}\mY^\infty_k(s)=0\,,\quad k\in[n-1]\,,\\
	&&\int_0^\infty \big(\mathbf{1}_{\{Q^\infty_k(s)>0\}}+\mathbf{1}_{\{Q^\infty_{k+1}(s)>0\}}\big)\,\mathrm{d}\modY^\infty_k(s)=0\,,
	\quad k\in[n-2]\,.
	\end{eqnarray}
		\end{enumerate}
\end{enumerate}
\end{prop}

\smallskip
The proof will follow  the ideas  in \cite{W1}, with additional complications caused by the last summand on the right-hand side of \eqref{skor_eq}. As there, for any $m\in\nn$ and any function $f\in D^m[0,\infty)$, we introduce the notation
\eq
\mathrm{Osc}_{t_1}^{t_2}(f)\,:=\,\sup_{t_1\leq s<t\leq t_2} \max_{1\leq k\leq m} |f_k(t)-f_k(s)|. 
\en  
The proof of Proposition \ref{inv_prin} is based on the following lemma.

\begin{lemma}\label{osc_lemma}
Suppose that the functions $\,q(\cdot),\,\centx (\cdot),\,\my (\cdot)\in D^{n-1}[0,\infty)$, $\mody (\cdot)\in D^{n-2}[0,\infty)$ fulfill the analogue 
$$
q(\cdot)=\centx(\cdot)+\mR\,\my (\cdot)+\modR\,\mody (\cdot)
$$
of \eqref{skor_eq}. 
% Yannis: I thought it made sense to phrase things this way. 
% (Yannis) with all capital letters replaced by small letters and all superscripts ``$N$'' removed. 
Moreover, suppose that the function $q (\cdot)$ takes values in $(\rr_+)^{n-1}$, $\my(0)=0$, $\mody(0)=0$, all components of $\my(\cdot) $, $\mody(\cdot)$ are non-decreasing and
\begin{eqnarray}
&&\int_0^\infty \mathbf{1}_{\{q_k(s)>0\}}\,\mathrm{d}\my_k(s)=0,\quad k\in[n-1], \label{loctime1}\\
&&\int_0^\infty \big(\mathbf{1}_{\{q_k(s)>0\}}+\mathbf{1}_{\{q_{k+1}(s)>0\}}\big)\,\mathrm{d}\mody_k(s)=0,\quad k\in[n-2], \label{loctime2}\\ 
&&\forall\,0\leq s<t:\;\mody_k(t)-\mody_k(s)\leq\min(\my_k(t)-\my_k(s),\my_{k+1}(t)-\my_{k+1}(s)).\label{loctime3} 
\end{eqnarray}  
Then, there is a constant $\,C\in (0, \infty)\,$ depending only on the entries of $\mR$ and $\modR$ $($but not on the particular functions $\,q (\cdot),\,\centx (\cdot), \,\my (\cdot), \,\mody(\cdot))\,$  such that
\eq\label{osc_ineq}
\forall\,~~0\leq t_1<t_2:\quad \mathrm{Osc}^{t_2}_{t_1}(q)+\mathrm{Osc}^{t_2}_{t_1}(\my)+\mathrm{Osc}^{t_2}_{t_1}(\mody)
\leq C\,\mathrm{Osc}^{t_2}_{t_1}(\centx).  
\en
\end{lemma}

\medskip

\noindent\textbf{Proof of Lemma \ref{osc_lemma}:} We proceed by induction over the dimension $n\geq2$. For $n=2$, the lemma is a direct consequence of Theorem 5.1 in \cite{W1}, since $\mody (\cdot)\equiv0$ in this case. From now on, we take $n\geq3$ and assume that the lemma holds for all $\nu=2,\ldots,n-1$. Consider first the case that there exists a $k\in[n-1]$ such that $\my_k(t_2)=\my_k(t_1)$. Then, \eqref{loctime3} shows that $\mody_k(t_2)=\mody_k(t_1)$ and $\mody_{k-1}(t_2)=\mody_{k-1}(t_1)$. Therefore, the induction hypothesis and the same argument as on pages 15-16 in \cite{W1} imply together that \eqref{osc_ineq} holds in this case.

\medskip

We now claim the existence of a vector $\lambda = (\lambda_1, \cdots, \lambda_n) \in(\rr_+)^{n-1}$ such that $(\lambda'\mR)_k\geq1$ and $(\lambda'\modR)_k\geq0$ for all $k\in[n-1]$, where $\lambda'$ stands for the transpose of $\lambda$. Indeed, noting that the off-diagonal elements in every column of $\mR$ are negative and add up to $-1$ (see Assumption \ref{coeff_asmp}(i)), and that the entries of each column of $\modR$ add up to $0$ and satisfy Assumption \ref{coeff_asmp}(ii), we see that we may choose $\lambda_1,\lambda_2,\ldots,\lambda_{n-1}$ as $f(a_1),f(a_2),\ldots,f(a_{n-1})$ for appropriate numbers $0=a_1<a_2<\ldots<a_{n-1}<1$ and an appropriate strictly concave function $f:\,[0,1]\rightarrow(0,\infty)$. Next, we fix such a $\lambda$ and proceed as in the derivation of the inequality (18) in \cite{W1} to conclude
\eq\label{y_bound}
\forall\,0\leq t_1<t_2:\quad \mathrm{Osc}_{t_1}^{t_2}(\my)
\leq\lambda'(q(t_2)-q(t_1))+\left(\sum_{k=1}^{n-1} \lambda_k\right)\mathrm{Osc}_{t_1}^{t_2}(\centx). 
\en 
Recalling \eqref{loctime3} and proceeding as in the derivation of the estimate (20) on page 17 of \cite{W1}, we can find a constant $C<\infty$ such that
\eq\label{mainestimate}
\forall\,0\leq t_1<t_2:\quad \mathrm{Osc}_{t_1}^{t_2}(q)
\leq C\Big(\mathrm{Osc}_{t_1}^{t_2}(\centx)+\max_{k\in[n-1]} q_k(t_2)\Big).
\en

\smallskip

The lemma can be now obtained from \eqref{mainestimate} in the same manner as Theorem 5.1 in \cite{W1} is obtained from the estimate (20) there. For the convenience of the reader, we give a short sketch of the proof. Fix $0\leq t_1<t_2$ and let $0<K<\infty$ be a constant, whose value will be determined later. We distinguish between two cases:
\begin{enumerate}[(a)]   
\item \quad $q_k(t_1)>K\mathrm{Osc}_{t_1}^{t_2}(\centx)$ for some $k\in[n-1]$.
\item \quad $q_k(t_1)\leq K\mathrm{Osc}_{t_1}^{t_2}(\centx)$ for all $k\in[n-1]$.
\end{enumerate}

\smallskip

In case (a), let $\tau$ be the first time $t\in[t_1,t_2]$ that $q_k(t)=0$ and set $\tau=\infty$ if the latter event does not occur. If $\tau=\infty$, then \eqref{osc_ineq} holds by the argument in the first paragraph of this proof as a consequence of the induction hypothesis. Now, suppose that $\tau\neq\infty$. Then, the same argument based on the induction hypothesis shows that
\eq
\mathrm{Osc}_{t_1}^{\tau}(q)\leq C\,\mathrm{Osc}_{t_1}^{\tau}(\centx), 
\en
where we have increased the value of $C<\infty$ if necessary. Next, let $C<K<\infty$. Since we are in case (a), we have
\eq
q_k(\tau)\geq q_k(t_1)-\mathrm{Osc}_{t_1}^{\tau}(q)\geq (K-C)\,\mathrm{Osc}_{t_1}^{\tau}(\centx)>0,
\en
which is a contradiction to $\tau\neq\infty$. 

\smallskip

In case (b), we distinguish two possibilities:
\begin{enumerate}[(i)]   
\item \quad $q_k(t)\leq K\mathrm{Osc}_{t_1}^{t_2}(\centx)$ for all $k\in[n-1]$ and $t\in[t_1,t_2]$.
\item \quad $q_k(t)>K\mathrm{Osc}_{t_1}^{t_2}(\centx)$ for some $k\in[n-1]$ and $t\in[t_1,t_2]$.
\end{enumerate}
In case (i), the inequality \eqref{osc_ineq} follows from 
\eq
\mathrm{Osc}_{t_1}^{t_2}(q)\leq \sup_{t_1\leq t\leq t_2} \max_{k\in[n-1]} q_k(t),
\en
\eqref{y_bound} and \eqref{loctime3}. In case (ii), we let $\tau$ be the first time $t\in[t_1,t_2]$ such that $q_k(t)>K\mathrm{Osc}_{t_1}^{t_2}(\centx)$ for some $k\in[n-1]$. Then, we split the time interval $[t_1,t_2]$ into $[t_1,\tau]$ and $[\tau,t_2]$, and argue as in case (i) on $[t_1,\tau]$ and as in case (a) on $[\tau,t_2]$. \ep

\bigskip

\noindent\textbf{Proof of Proposition \ref{inv_prin}:} First, we claim that part (a) of Proposition \ref{inv_prin} is a consequence of the inequality \eqref{osc_ineq} in Lemma \ref{osc_lemma}. Indeed, the necessary and sufficient conditions (a) and (b) of Corollary 3.7.4 in \cite{EK} hold for any subsequence of $\centX^N(\cdot)$, $N>0$, and carry over to the same subsequence of $(Q^N(\cdot),\centX^N(\cdot),\mY^N(\cdot),\modY^N(\cdot))$, $N>0$ via the inequality \eqref{osc_ineq}. For more details on the same argument, please see the proof of Theorem 4.1 in \cite{W1}. 

\smallskip
Now, we turn to the proof of part (b) of Proposition \ref{inv_prin} and let $\,(Q^\infty (\cdot),\centX^\infty (\cdot),\mY^\infty (\cdot),$ $ \modY^\infty(\cdot))$ be a limit point as there. The properties (i) and (ii) for it can be seen by using the Skorokhod Representation Theorem in the form of Theorem 3.1.8 in \cite{EK} (noting the path continuity of the limit) for the subsequence of $(Q^N(\cdot),\centX^N(\cdot),\mY^N(\cdot),\modY^N(\cdot))$, $N>0$, which converges to 
that limit point, and taking the almost sure limit on both sides of the identity \eqref{skor_eq} through this particular subsequence. To deduce property (iii), we let $g:\,[0,\infty)\rightarrow[0,1]$ be a continuous function such that, for some $\delta>0$, it holds $g(a)=0$ whenever $0\leq a\leq\delta$, and $g(a)=1$ whenever $a\geq2\delta$. By arguing as in the second half of the proof of Theorem 4.1 in \cite{W1}, we conclude that the quantities
\begin{eqnarray}
&&\int_0^\infty g\left(Q^N_k(s)\right)\,\mathrm{d}\mY^N_k(s)=0,\quad k\in[n-1], \\
&&\int_0^\infty g\left(Q^N_k(s)\right)+g\left(Q^N_{k+1}(s)\right)\,\mathrm{d}\modY^N_k(s)=0,\quad k\in[n-2]
\end{eqnarray} 
converge to the quantities
\begin{eqnarray}
&&\int_0^\infty g\left(Q^\infty_k(s)\right)\,\mathrm{d}\mY^\infty_k(s),\quad k\in[n-1], \label{ltlim1}\\
&&\int_0^\infty g\left(Q^\infty_k(s)\right)+g\left(Q^\infty_{k+1}(s)\right)\,\mathrm{d}\modY^\infty_k(s),\quad k\in[n-2] \label{ltlim2}
\end{eqnarray} 
when we pass to the limit through the same subsequence as before. Now, letting $\, \delta \downarrow 0\,$, we obtain the property (iii). \ep  

\smallskip

We are now ready for the proof of Theorem \ref{thm_conv}. 

\medskip

\noindent\textbf{Proof of Theorem \ref{thm_conv}: Step 1.} Consider the family of processes $\centX^N(\cdot)$, $N>0$ in \eqref{eq:whatisxN}. The family of processes $\Big(\frac{1}{\sqrt{N}}\,\centS_k^L(Nt),\,t\geq0,\,\frac{1}{\sqrt{N}}\,\centS_k^R(Nt),\,t\geq0,\,k\in[n]\Big)$, $N>0$, without the time change, is tight by Lemma \ref{lem:BMconv}. Moreover, the time-changes in \eqref{eq:whatisxN} are Lipschitz functions of time. That is, there exists a constant $\Theta<\infty$ such that 
\[
\forall\,0\leq s<t:\quad \max\left(T_k^L(t)-T_k^L(s),T_k^R(t)-T_k^R(s)\right)\leq \Theta(t-s),\quad k\in[n].
\]
Hence, using the necessary and sufficient conditions for tightness of Corollary 3.7.4 in \cite{EK}, the tightness of $\centX^N(\cdot)$, $N>0$
is easily verified. However, we still have to identify the limit points.   
% Yannis: I thought this phraseology was more to the point. Please verify.

\medskip

\noindent\textbf{Step 2.} At this stage, we can use Proposition \ref{inv_prin} to conclude that the family $\left(Q^N(\cdot),\centX^N(\cdot),\mY^N(\cdot),\modY^N(\cdot)\right)$, $N>0$ is tight. Now, recall that, for any $N>0$, $\mY^N(\cdot)$ and $\modY^N(\cdot)$ can be expressed, as in \eqref{eq:whatisy}, in terms of the times that the particles spend in collisions. The tightness of $\,\mY^N(\cdot)$, $N>0\,$ and $\,\modY^N(\cdot)$, $N>0\,$ now shows that the processes 
\eq\label{eq:loctimenegligible}
\frac{1}{N} I_k\left(Nt\right),\;t\geq0,\quad k\in[n-1],\quad \frac{1}{N} I_{k,k+1}(Nt),\;t\geq0,\quad k\in[n-2]
\en
all tend to zero in $D[0,\infty)$.

\medskip

\noindent\textbf{Step 3.} Putting together the conclusion of Step 2, \eqref{eq:whatisb} and \eqref{eq:whatisb2}, we conclude that each of the processes
\[
\frac{1}{N} T_k^L\left(Nt\right),\;t\geq0,\quad k\in[n],\quad \frac{1}{N} T_k^R\left(Nt\right),\;t\geq0,\quad k\in[n] 
\] 
converge  to the process $t$, $t\geq0$ in $D[0,\infty)$. With the help of the Lemma on page 151 in \cite{B99}, preceding Theorem 14.4,  we deduce f  that the joint distributions of 
\[
\left(\frac{1}{\sqrt{N}}\,\centS_k^L(T_k^L(Nt)),\,t\geq0,\;\frac{1}{\sqrt{N}}\,\centS_k^R(T_k^R(Nt)),\,t\geq0,\; k\in[n]\right), 
\]
now  {\sl with the time-change,} converge on $D^{2n}[0,\infty)$ to the limiting distribution  described in Lemma \ref{lem:BMconv} (note that we can improve the topology used in \cite{B99} to the topology of uniform convergence on compacts by observing the path continuity of the limit process, modifying the paths of the jump processes to continuous, piecewise linear paths and using the fact that the Skorokhod topology relativized to the space of continuous functions coincides with the locally uniform topology there).  

\medskip

\noindent\textbf{Step 4.} Step 3 implies that, in the limit $N\rightarrow\infty$, the processes $\centX^N(\cdot)$, $N>0$ converge vaguely in $D^{n-1}[0,\infty)$ to a multidimensional Brownian motion with drift and diffusion co\"efficients as described in Theorem \ref{thm_conv}. To conclude the proof, we note that every limit point $\,(Q^\infty (\cdot),\centX^\infty (\cdot),\mY^\infty (\cdot),\modY^\infty(\cdot))\,$ of the collection $\,(Q^N(\cdot),\centX^N(\cdot),\mY^N(\cdot),$  $\modY^N(\cdot))\,$, $\,N>0\,$ has continuous paths by \eqref{osc_ineq} and the fact that the processes $\centX^N(\cdot)$, $N>0$ converge to a process with continuous paths. Thus, by part (b) of Proposition \ref{inv_prin}, every such limit point  satisfies the properties (i), (ii), (iii) there, with $\centX^\infty (\cdot)$ being the multidimensional Brownian motion just described. Lastly, by Lemma \ref{genRW}, we can identify $Q^\infty (\cdot)$ with $\xi (\cdot)$, $\mY^\infty (\cdot)$ with the boundary local times of $\xi (\cdot)$, and deduce that $\modY^\infty (\cdot) \equiv 0$. This completes the proof. \ep

\bigskip

Finally, we consider the limit of the entire collection of jump processes   
% (Yannis) random walks
 $\left(\Gamma_k(\cdot),k\in[n]\right)$ (we refer the reader to the expression \eqref{eq:existence}). Let $\Gamma^N(\cdot)$ denote the vector of centered and scaled jump processes   
% (Yannis)   walks 
given, for every $k\in[n]$, by
\[
\Gamma_k^N(t):=\frac{1}{\sqrt{N}}\Gamma_k(0)+\frac{1}{\sqrt{N}}\left[S^R_k\left(T^R_k(Nt)\right)-S^L_k\left(T^L_k(Nt)\right)\right]-(a-b)t\sqrt{N}, \;\;\;t\geq0\,.
\]

\begin{thm}\label{thm:conv_full}
Suppose that Assumption \ref{coeff_asmp} holds and that $\lim_{N\rightarrow\infty} N^{-1/2}\Gamma_k(0)=R(0)\in\mathbb{W}^n$ in distribution. Further, let $R(\cdot)=\left(R_1(\cdot),R_2(\cdot),\ldots,R_n(\cdot)\right)$ denote the continuous $n$-dimensional semimartingale taking values in $\mathbb{W}^n$ and satisfying 
% Yannis  
%\footnote{~In the expression below, as well as in (\ref{disc_ranks}), I have changed the last subscript in the last denominator from $k+1$ to $k$. This way condition (\ref{elastic}) is satisfied; we may wish to call the reader's attention to this fact, will leave this to your judgment. (IK)}
\begin{eqnarray*}
R_k(t)&=&R_k(0)+ \left(\lambda_k^R-\lambda_k^L \right)t+\left(a^3\left(\sigma^R_k\right)^2+b^3\left(\sigma^L_k\right)^2\right)^{1/2}\betab_k(t) \\
&& +\frac{a\left(\theta_k^R-1\right)+b}{a\theta_k^R+b\theta_{k+1}^L}\Lambda^{(k,k+1)}(t)
-\frac{b\left(\theta_k^L -1\right)+a}{a\theta^R_{k-1}+b\theta^L_{k }}\Lambda^{(k-1,k)}(t),\quad t\geq0,
\end{eqnarray*}
$k\in[n]$ with the same notation as in \eqref{ranks}. Then, the processes $\Gamma^N(\cdot)$, $N>0$ converge in $D^n[0,\infty)$ to the process $R(\cdot)$ described above in the limit $N\rightarrow\infty$.  
\end{thm}

\medskip

\noindent\textbf{Proof.} The main observation is that, for any fixed $N>0$, we have
\eq\label{disc_ranks}
\begin{split}
\Gamma_k^N(t)=&\,\frac{1}{\sqrt{N}}\Gamma_k(0)+\frac{1}{\sqrt{N}}\left[\centS^R_k\left(T^R_k(Nt)\right)
-\centS^L_k\left(T^L_k(Nt)\right)\right]\\
&+\frac{a\left(\theta_k^R-1\right)+b}{a\theta_k^R+b\theta_{k+1}^L}\,\mY^N_k(t)
-\frac{b\left(\theta_k^L -1\right)+a}{a\theta^R_{k-1}+b\theta^L_{k }}\,\mY^N_{k-1}(t)\\
&+\left(-a\left(\theta^R_k-1\right)+b\left(\theta^L_k-1\right)\right)\,\modY^N_{k-1}(t),\quad t\geq0,
\end{split}
\en
$k\in[n]$. The same steps as in the proof of Theorem \ref{thm_conv} now show that the processes in the first line of \eqref{disc_ranks} converge  jointly, in distribution, to the components of a multidimensional Brownian motion with drift and diffusion co\"efficients as in the statement of this theorem, the process $\mY^N(\cdot)$ converges to the process of boundary local times of a reflected Brownian motion as in Theorem \ref{thm_conv}, whereas the process $\modY^N(\cdot)$ converges to zero. Therefore, the processes $\Gamma^N(\cdot)$, $N>0$ must also converge in distribution, and one can identify the limit of the ``noise part'' with the appropriate multidimensional Brownian motion and the limit of the ``local time part'' with the local time part in the decomposition of the process $R(\cdot)$ in the statement of the theorem. This completes the argument. \ep 

%%%%%%%%%%%%%%%%%%%%%%%%%%%%%%%%%%%%%%%
\section{Additional Determinantal Structures}
%%%%%%%%%%%%%%%%%%%%%%%%%%%%%%%%%%%%%%%

This last section studies conditions on the parameters $\,b_1,b_2,\ldots,b_n\,$, $\sigma_1,\sigma_2,\ldots,\sigma_n$ and $\,q^{\pm}_1,q^{\pm}_2,\ldots,q^{\pm}_n\,$, under which the process $\,R(\cdot)=(R_1(\cdot),R_2(\cdot),\ldots,R_n(\cdot))\,$ of ranks as in \eqref{ranks} has a probabilistic structure of {\it determinantal type,} in the sense that its transition densities are of the generalized \textsc{Karlin-McGregor} form
\eq\label{det_density}
p(t,r,\tilde{r})\,=  \sum_{\sigmab \in \mathbb{S}_n} \kappa_\sigmab\,\prod_{k=1}^n \, f^{\,k,\sigmab(k)} \big(t,\tilde{r}_{\sigmab(k)}-r_k\big)\,.
\en 
Here $\,\mathbb{S}_n\,$ is the group of permutations of a set with $n$ elements, whereas $\,\kappa_\sigmab\,$, $\sigmab \in \mathbb{S}_n$ and $f^{k,\ell}$, $1\leq k,\,\ell \leq n$ are suitable real numbers and real-valued functions, respectively. 

\smallskip

This question is motivated by the following two extreme cases: 

\smallskip
\noindent
{\it (i)}  If one considers $\,b_1=b_2=\ldots=b_n\,$, $\,\sigma_1=\sigma_2=\ldots=\sigma_n\,$ and $q^{\pm}_1=q^{\pm}_2=\ldots=q^{\pm}_n=1/2$, then the ranks evolve as the ordered system of $n$ Brownian motions, each with drift $b_1$ and dispersion  $\sigma_1$. Therefore, in this case
\eq
p(t,r,\tilde{r})\,=\sum_{\sigmab \in S_n} \, \prod_{k=1}^n \, \varphi_{\,b_1,\sigma_1 }\big(t,\tilde{r}_{\sigmab(k)}-r_k\big)\,,
\en 
where $\varphi_{\, b_1,\sigma_1 }(t,\cdot)$ denotes the Gaussian density with mean $\,b_1\,t\,$ and variance $\,\sigma_1^2\,t\,$. 

\smallskip
\noindent
{\it (ii)} Now, consider the case of $\,b_1=b_2=\ldots=b_n=0\,$, $\,\sigma_1=\sigma_2=\ldots=\sigma_n=1\,$ and 
\eq
q_k^+=1,\quad k=1,2,\ldots,n\,,\quad q_k^-=0,\quad k=1,2,\ldots,n\,.
\en
Then, the process of ranks is given by the continuous version of the totally asymmetric simple exclusion process (TASEP) treated in detail by \textsc{Warren} (2007) in section 4 of \cite{Wa}. In this case, the transition probability densities are of the form
\eq
p(t,r,\tilde{r})=\sum_{\sigmab \in \mathbb{S}_n} (-1)^{\mathrm{sgn}(\sigmab)} \prod_{k=1}^n \, \psi^{k,\sigmab (k)} \big(t,\tilde{r}_{\sigmab(k)}-r_k\big)
\en
for suitable functions $\,\psi^{k,\ell}\,$, $1\leq k,\,\ell \leq n$ (see Proposition 8 in \cite{Wa}) and with $\mathrm{sgn}(\sigmab)$ standing for the signum of a permutation $\, \sigmab \in \mathbb{S}_n\,$. 

\smallskip
The main result of this section is the following proposition.

\begin{prop}
Suppose that, for every $\epsilon>0$, the transition probability densities of the process in \eqref{ranks} belong to the function space $C_b((\epsilon,\infty)\times\mathbb{W}^n\times\mathbb{W}^n)$, and are continuously differentiable in the first co\"ordinate and twice continuously differentiable in the second co\"ordinate with derivatives in $C_b((\epsilon,\infty)\times\mathbb{W}^n\times\mathbb{W}^n)$. 

\smallskip
Then these transition densities are given by \eqref{det_density} with suitable real constants $\kappa_\sigmab\,$, $\sigmab\in \mathbb{S}_n\,$ and real-valued functions $\,f_{k,\ell}\,$, $\,1\leq k,\,\ell \leq n\,$, if and only if:  $\,b_1=b_2=\ldots=b_n$, $\sigma_1=\sigma_2=\ldots=\sigma_n\,$ and for each $\,k=1,2,\ldots,n-1\,$  we have, either  

\smallskip
\noindent
 (i)~~ $\,\,\,q_k^-=q_{k+1}^+=1/2\,$, or \\ (ii) ~~$\,\,q_k^-=0$, $q_{k+1}^+=1\,$, or \\ (iii) ~~ $\,q_k^-=1$, $q_{k+1}^+=0\,$. 

\medskip
\noindent
Moreover, if this is the case, one may choose $\,\kappa_\sigmab=1\,$, $\sigmab\in \mathbb{S}_n\,$ and the functions $f^{k,\ell}$, $1\leq k,\,\ell\leq n$ according to the formulas \eqref{form1}-\eqref{form5} below.  
\end{prop}

\medskip
\noindent{\bf Proof.} We pick a bounded and continuous function $\,g:\,\mathbb{W}^n\rightarrow\rr\,$ vanishing in a neighborhood of the boundary of $\,\mathbb{W}^n\,$,   an arbitrary real number $\, T >0\,$, and consider the martingale
\eq\label{Rmart}
\int_{\mathbb{W}^n} p(T-t,R(t),\tilde{r})\,g(\tilde{r})\,\mathrm{d}\tilde{r}=:F(T-t,R(t)),\quad  0 \le t \le T\,.  
\en
Since the change of variables formula of \cite{HR} is applicable to the spacings process $Z(\cdot)$, the construction of the process $\,R(\cdot)\,$ in section 2 shows that an analogous change of variables formula holds for the process $\,R (\cdot)\,$. Applying this formula to the right-hand side of \eqref{Rmart} and differentiating under the integral, we obtain from the martingale property of the process in \eqref{Rmart} the heat equation 
\begin{equation} 
  \partial_{t} \, p (t,r,\tilde{r})\,=\,\frac{1}{\,2\,}\sum_{k=1}^n \sigma_k^2\, \partial_{r_k}^2 \,p (t,r,\tilde{r})
  +\sum_{k=1}^n b_k\, \partial_{r_k} \,p (t,r,\tilde{r})   \,, \label{heat_eq}
\end{equation}
and the elastic boundary condition 
\begin{equation} 
  q^-_k\, \partial_{r_k}\,p (t,r,\tilde{r})- q^+_{k+1}\,\partial_{r_{k+1}} \,p (t,r,\tilde{r})=0\quad\mathrm{whenever}\;r_k=r_{k+1}\,, 
\label{bnd_cond}
\end{equation}

\medskip
\noindent
for $\, t>0\,,\;(r,\tilde{r})\in (\mathbb{W}^n)^2\,$.

\medskip

Substituting the expression of \eqref{det_density} into \eqref{bnd_cond}, we deduce  
\begin{eqnarray*}
&&q_k^- \sum_{\sigmab\in \mathbb{S}_n} \kappa_\sigmab \,\prod_{\ell\neq k} f^{\, \ell,\, \sigmab (\ell)}(t,\tilde{r}_{\sigmab (\ell)}-r_\ell)\,
D_2 f ^{\,k,\,\sigmab (k)}(t,\tilde{r}_{\sigmab(k)}-r_k)\\
&-&q_{k+1}^+ \sum_{\sigmab \in \mathbb{S}_n} \kappa_\sigmab \,
\prod_{\ell\neq k+1} f^{\, \ell,\,\sigmab (\ell)}(t,\tilde{r}_{\sigmab (\ell)}-r_\ell)\,
D_2 f^{\,k+1,\,\sigmab(k+1)}(t,\tilde{r}_{\sigmab(k+1)}-r_{k+1})\,=\,0
\end{eqnarray*}
whenever $\,r_k=r_{k+1}\,$. Here, $\,D_2\,$ denotes   differentiation with respect to the second argument. Plugging in $r_k$ for $r_{k+1}$ in this last expression, and grouping together functions that have the same arguments, results in  
\eq\label{bdry1}
\begin{split}
\forall ~~1\leq k\leq n,\,1\leq\ell_1\neq \ell_2\leq n,
&\,\sigmab(k)=\ell_1,\,\sigmab(k+1)=\ell_2,\,\widetilde{\sigmab}(k)=\ell_2,\,\widetilde{\sigmab}(k+1)=\ell_1:\\
0\,=&\,\,q_k^-\kappa_\sigmab\,D_2 f^{k,\ell_1} (t,\widetilde{r}_{\ell_1}-r_k) \,f^{k+1,\ell_2}(t,\widetilde{r}_{\ell_2}-r_k)\\
&+q_k^-\kappa_{\widetilde{\sigmab}} f^{k+1,\ell_1}(t,\widetilde{r}_{\ell_1}-r_k) \,D_2 f^{k,\ell_2}(t,\widetilde{r}_{\ell_2}-r_k)\\
&-q_{k+1}^+\kappa_\sigmab\,f^{k,\ell_1}(t,\widetilde{r}_{\ell_1}-r_k) D_2 f^{k+1,\ell_2}(t,\widetilde{r}_{\ell_2}-r_k)\\
&-q_{k+1}^+\kappa_{\widetilde{\sigmab}} D_2 f^{k+1,\ell_1}(t,\widetilde{r}_{\ell_1}-r_k) \,f^{k,\ell_2}(t,\widetilde{r}_{\ell_2}-r_k). 
\end{split}
\en

\medskip
\noindent
Let us recall now \eqref{elastic}, and take the Fourier transform with respect to the variables $\tilde{r}_{\ell_1}-r_k$ (parameter $a$) and $\tilde{r}_{\ell_2}-r_k$ (parameter $b$), to obtain equations of the form
\eq
\begin{split}
q\,\big(a\,G(t,a)\,H(t,b)+K(t,a)\,b\,L(t,b)\big) \quad\quad\quad\quad\quad\quad\\  
=\,(1-q)\big(G(t,a)\,b\,H(t,b)+a\,K(t,a)\,L(t,b)\big) \,,
\end{split}
\en
or equivalently
\eq
G(t,a)\,H(t,b)\,\frac{\,\,q\,a+(1-q)\,b\,\,}{q\,b+(1-q)\,a}\,=\,K(t,a)\,L(t,b)\,. 
\en
The only cases in which the fraction $$\frac{q\,a+(1-q)\,b}{q\,b+(1-q)\,a}$$ can be factored as the product of a function only of $a$ and a function only of $b$, are given by $q=1/2$, $q=1$ and $q=0$. Moreover, up to multiplicative constants, 
\begin{eqnarray}
&&f^{k,\ell}=f^{k+1,\ell}\,,\quad \ell=1,2,\dots,n\quad\mathrm{whenever}\;q_k^-=1/2\,, \label{bndr_cond1}\\
&&D_2 f^{k,\ell}=f^{k+1,\ell}\,,\quad \ell=1,2,\dots,n\quad\mathrm{whenever}\;q_k^-=1\,, \label{bndr_cond2}\\
&&f^{k,\ell}=D_2 f^{k+1,\ell}\,,\quad \ell=1,2,\dots,n\quad\mathrm{whenever}\;q_k^-=0\,. \label{bndr_cond3} 
\end{eqnarray}

\medskip

Next, we introduce the linear parabolic operators
\begin{eqnarray}
\mathcal R_k\,=\,\partial_t-\frac{1}{\,2\,}\,\sigma_k^2\,\partial_{r_k}^2-b_k\,\partial_{r_k}\,,\quad k=1,2,\ldots,n  
\end{eqnarray}
and substitute the expression of \eqref{det_density} into \eqref{heat_eq}, to obtain  
\eq
0\,=\sum_{\sigmab \in \mathbb{S}_n} \kappa_\sigmab\,\sum_{k=1}^n \Big(\prod_{\ell\neq k} \, f^{\,\ell,\sigmab(\ell)} \big(t,\tilde{r}_{\sigmab(\ell)}-r_\ell \big)\Big)\,
\,\mathcal R_k\,f^{\,k,\sigmab(k)}(t,\tilde{r}_{\sigmab(k)}-r_k)\,,
\en
which shows 
\eq
\mathcal R_k\,f^{\,k,\sigmab(k)}=0,\quad k=1,2,\ldots,n\,,\;\;\sigmab\in\mathbb{S}_n\,. 
\en
This and \eqref{bndr_cond1}-\eqref{bndr_cond3} imply $b_1=b_2=\ldots=b_n$ and $\sigma_1=\sigma_2=\ldots=\sigma_n$. This completes the proof of the ``only if'' part.

\medskip

Conversely, if the conditions on the drift, dispersion and collision parameters in the proposition are satisfied, we define $f^{k,\ell}$, $1\leq k,\,\ell\leq n$ by 
\begin{eqnarray}
&&f^{k,k}=\varphi_{b_1,\sigma_1}\,,\quad k=1,2,\ldots,n\,, \label{form1}\\
&&f^{k,\ell}=f^{k+1,\ell}\,,\quad \ell=1,2,\ldots,n\quad\mathrm{whenever}\;\,q_k^-=1/2\,,\;k=1,2,\ldots,n-1\,,\\
&&D_2f^{k,\ell}=f^{k+1,\ell}\,,\quad \ell=1,2,\ldots,n\quad\mathrm{whenever}\;\,q_k^-=1,\;k=1,2,\ldots,n-1\,,\\
&&f^{k,\ell}=D_2f^{k+1,\ell}\,,\quad \ell=1,2,\ldots\,,n\quad\mathrm{whenever}\;\,q_k^-=0,\;k=1,2,\ldots,n-1\,. \label{form5}
\end{eqnarray}
We can express this state of affairs as follows: In order to  determine the entry $\,f^{k,\ell}(t, \cdot\,)\,$ in \eqref{det_density} for $\ell>k$, we count the number $\,\mathfrak{u}\,$ of ones in $\{q_k^-,\ldots,q_{\ell-1}^-\}$ and the number $\,\mathfrak{z}\,$ of zeros in $\{q_k^-,\ldots,q_{\ell-1}^-\}$; and then compute $\,f^{k,\ell}(t, \cdot\,)\,$  by differentiating the Gaussian probability density function  $\,\varphi_{b_1,\sigma_1} (t, \cdot)\,$ with respect to its second co\"ordinate $\,\mathfrak{u}\,$ times, and integrating the result with respect to the second co\"ordinate $\,\mathfrak{z}\,$ times. The entries $\,f^{k,\ell}(t, \cdot\,)\,$ for $\ell<k$ are computed similarly. 

Now, in the case 
\eq
q_k^-\,=\,1/2\,,\quad k=1,2,\ldots,n\,,
\en
one just needs to argue as in the beginning of this section to finish the proof; whereas in all other cases one can argue as in the proofs of Proposition 8 and Lemma 7 in \cite{Wa} to deduce that the expression of \eqref{det_density} with these choices of $\,f^{k,\ell}(t, \cdot)\,$, $1\leq k,\,\ell\leq n$ and $\kappa_\sigmab=1$, $\sigmab\in \mathbb{S}_n$ gives the transition densities of the process $\,R(\cdot)\,$. This completes the proof. \ep

\bigskip

\bibliographystyle{alpha}

\bigskip

\end{document}